\documentclass[11pt,french]{amsart}
\usepackage[francais]{babel}
\usepackage[utf8]{inputenc}
\usepackage{latexsym}
\usepackage[all]{xy}
\usepackage{amsmath}
\usepackage{amssymb}
\usepackage{amsfonts}

\usepackage{color}
\definecolor{niceblue}{rgb}{0,0,0.6}
\usepackage{latexsym}

\usepackage{stmaryrd}
\usepackage{amsmath, amscd, amssymb, mathrsfs, accents, amsfonts}
\usepackage{url}
\usepackage[all]{xy}
\usepackage{longtable}
\usepackage{dsfont}
\usepackage{tikz}
\usetikzlibrary{decorations.pathmorphing}
\usetikzlibrary{calc}
\usetikzlibrary{decorations.markings}

\def\N{{\mathbb N}}
\def\Z{{\mathbb Z}}
\def\Q{{\mathbb Q}}
\def\R{{\mathbb R}}
  
\def\cP{{\mathcal P}}
\def\cA{{\mathcal A}}
\def\cC{{\mathcal C}}

\def\cS{{\mathcal S}}

\def\equiv#1{\smash{\mathop{\sim}\limits_{#1}}}

\def\tf{$\sqcup\hskip-2.3mm\sqcap$}   
\def\llrightarrow{{\relbar\joinrel\relbar\joinrel\rightarrow}} 
{\newtheorem{theorem}{Th\'eor\`eme}[section]   
   \newtheorem{proposition}[theorem]{Proposition}
   \newtheorem{lemma}[theorem]{Lemme}  
   \newtheorem{remark}[theorem]{Remarque} 
   \newtheorem{corollary}[theorem]{Corollaire}

   \newtheorem*{definition}{D\'efinition}    
{\theoremstyle{definition}   
   }
{\theoremstyle{remark}

   \newtheorem*{preuve}{Preuve}   
        
   \newtheorem*{preuveth}{Preuve du th\'eor\`eme}    
     
\author{Michel Vaqui\'e}
\address{Institut de Math\'ematiques de Toulouse UMR 5219, CNRS, Universit\'e de Toulouse,  
UPS, 118 route de Narbonne, F-31062 Toulouse Cedex 9, 
France} 
\email{vaquie@math.univ-toulouse.fr} 
\title[Valuation augment\'ee et paire minimale]{Valuation augment\'ee et paire minimale} 

   \thanks{$*$ Partially supported by the grant of the Agence Nationale de la Recherche “CatAG”ANR-17-CE40-0014.} 

\begin{document}

\begin{abstract} 
Soit $(K,\nu )$ un corps valu\'e, les notions de \emph{valuation augment\'ee}, de \emph{valuation augment\'ee limite} et de \emph{famille admise} de valuations permettent de donner une description de toute valuation $\mu$ de $K[x]$ prolongeant $\nu$.  
Dans le cas o\`u le corps $K$ est alg\'ebriquement clos cette description est particuli\`erement simple et nous pouvons la r\'eduire aux notions de \emph{paire minimale} et de \emph{famille pseudo-convergente}. 

Soient $(K,\nu )$ un corps valu\'e hens\'elien et $\bar\nu$ l'unique extension de $\nu$ \`a la cl\^oture alg\'ebrique $\bar K$ de $K$ et soit $\mu$ une valuation de $K[x]$ prolongeant $\nu$, nous \'etudions les extensions $\bar\mu$ de $\mu$ \`a $\bar K[x]$ et nous donnons une description des valuations $\bar \mu _i$ de $\bar K[x]$ qui sont les extensions des valuations $\mu _i$ appartenant \`a la famille admise associ\'ee \`a $\mu$.  

\vskip .2cm

\noindent {\scshape Abstract}.  Let $(K, \nu)$ be a valued field, the notions of \emph{augmented valuation}, of \emph{limit augmented valuation} and of \emph{admissible family} of valuations enable to give a description of any valuation $\mu$  of $K [x]$ extending $\nu$. In the case where the field $K$ is algebraically closed, this description is particularly simple and we can reduce it to the notions of \emph{minimal pair} and \emph{pseudo-convergent family}.

Let $(K, \nu )$ be a henselian valued field and $\bar\nu$ the unique extension of $\nu$ to the algebraic closure $\bar K$ of $K$ and let $\mu$ be a valuation of $ K [x]$ extending $\nu$, we study the extensions $\bar\mu$ from $\mu$ to $\bar K [x]$ and we give a description of the valuations $\bar\mu _i$ of $\bar K [x]$ which are the extensions of the valuations $\mu _i$ belonging to the admissible family associated with $\mu$.
 \end{abstract}


\vskip .2cm 
\date{Mai 2020}

\maketitle   

\tableofcontents

		      \section*{Introduction}
%

Soit $K$ un corps muni d'une valuation $\nu$, nous pouvons obtenir toute valuation ou pseudo-valuation $\mu$ de l'anneau des polyn\^omes $K[x]$ qui prolonge $\nu$ gr\^ace \`a {\it une famille admise de valuations} ${\mathcal A}=\bigl ( \mu _i \bigr )_{i \in I}$, o\`u l'ensemble $I$ est un ensemble totalement ordonn\'e (cf. th\'eor\`emes 2.4 et 2.5 de \cite{Va 1}). 
De plus chaque valuation $\mu _i$ de la famille est obtenue comme valuation augment\'ee ou comme valuation augment\'ee limite associ\'ee \`a un polyn\^ome-cl\'e ou un polyn\^ome-cl\'e limite $\phi _i$. 

La famille ${\mathcal A}=\bigl ( \mu _i \bigr )_{i \in I}$ converge vers la valuation $\mu$ dans le sens o\`u pour tout polyn\^ome $f$ de $K[x]$ la famille de valeurs $\bigl ( \mu _i  (f) \bigr )_{i \in I}$ est croissante et v\'erifie $\mu (f) = Sup \bigl ( \mu _i  (f) \ ; \ i \in I\bigr )$. 
En particulier si la famille $I$ a un plus grand \'el\'ement  $\bar\iota$ la valuation $\mu$ est \'egale \`a la valuation $\mu _{\bar\iota}$ et il existe un polyn\^ome $\phi = \phi _{\bar\iota}$ qui d\'efinit la valuation $\mu$ comme valuation augment\'ee ou valuation augment\'ee limite. 
Nous disons dans ce cas que la valuation est \emph{bien sp\'ecifi\'ee} et que le polyn\^ome $\phi$ \emph{d\'efinit} la valuation. 
Par d\'efinition ce polyn\^ome appara\^\i t dans la construction de la famille admise ${\mathcal A}=\bigl ( \mu _i \bigr )_{i \in I}$ associ\'ee \`a la valuation $\mu$. 

\vskip .2cm 

Comme tout polyn\^ome-cl\'e ou polyn\^ome-cl\'e limite est un polyn\^ome irr\'eductible de $K[x]$ dans le cas o\`u le corps $K$ est alg\'ebriquement clos les seuls polyn\^omes-cl\'es sont de degr\'e un et les familles admises sont particuli\`erement simples: si la valuation $\mu$ est bien sp\'ecifi\'ee elle est d\'efinie par un polyn\^ome-cl\'e $\phi$ de la forme $\phi (x) = x-a$, sinon elle est d\'efinie par une famille infinie de polyn\^omes $(\phi _{\alpha})$ de la forme $\phi _{\alpha}(x) = x-a _{\alpha}$. 
Plus pr\'ecis\'ement, comme $K$ est alg\'ebriquement clos la valuation $\mu$ est enti\`erement d\'etermin\'ee par les valeurs prises pour les polyn\^omes $f$ de la forme $f(x)=(x-b)$, et nous avons dans le cas o\`u $\mu$ est bien sp\'ecifi\'ee  
$\mu (x-b) \ = \ Inf \left ( \nu (a -b)  , \delta \right )$ et la valuation $\mu$ est associ\'ee \`a une \emph{paire minimale}, 
et dans le cas d'une famille infinie la valuation $\mu$ est la valuation associ\'ee \`a la \emph{famille pseudo-convergente} $(a _{\alpha })$ ({\bf Proposition \ref{prop:corpsalgclos}}).   
Il est aussi possible de d\'ecrire une valuation $\mu$ de $K[x]$ par une boule ferm\'ee ou par une famille d\'ecroissante de boules ferm\'ees dans $K$ pour la distance ultram\'etrique associ\'ee \`a la valuation $\nu$ de $K$  ({\bf Proposition \ref{prop:boules}}).  

\vskip .2cm 

Soient $(K,\nu )$ un corps valu\'e quelconque et $\mu$ une valuation de $K[x]$ prolongeant $\nu$, alors si $\bar\nu$ est une extension de $\nu$ \`a la cl\^oture alg\'ebrique $\bar K$ de $K$ il existe une extension $\bar\mu$ de $\mu$ \`a $\bar K[x]$ qui prolonge $\bar\nu$. 
Soit $\mathcal{A} = \bigl ( \mu _i \bigr ) _{i \in I}$ la famille admise associ\'ee \`a $\mu$, nous voulons d\'ecrire les valuations $\bar\mu _i$ de $\bar K[x]$ obtenues comme prolongement des valuations $\mu _i$ appartenant \`a la famille $\mathcal A$, et les familles de boules ferm\'ees de $\bar K$ associ\'ees aux valuations $\bar\mu _i$. 

Dans le cas o\`u $(K,\nu )$ est un corps valu\'e hens\'elien nous associons \`a la famille admise associ\'ee \`a la valuation $\mu$ une famille d\'ecroissante $\bigl ( {\mathcal B} _i \bigr )_{i \in I}$ de r\'eunions finies de boules ferm\'ees de $\bar K$. 
Plus pr\'ecis\'ement chaque ${\mathcal B} _i $ est la r\'eunion des boules disjointes $B_{i}^{(r)}$, le groupe de Galois agit transitivement sur l'ensemble fini $\{ B_{i}^{(r)} \}$, et  pour tout $i<j$ dans $I$ chaque boule $B_{i}^{(r)}$ de ${\mathcal B} _i$ contient $s$ boules $B_j^{(l)}$ appartenant \`a ${\mathcal B} _j$, o\`u $s$ est un entier ind\'ependant de la boule $B_i^{(r)}$ choisie, et toute boule  $B_j^{(l)}$ de ${\mathcal B} _j$ est contenue dans une boule $B_i^{(r)}$ de ${\mathcal B} _i$ ({\bf Proposition \ref{prop:familles-incluses}}).  

L'intersection des ${\mathcal B} _i$ est un sous-ensemble ${\mathbf B}(\mu )$ de $\bar K$, appel\'e \emph{ensemble caract\'eristique de la valuation} $\mu$, cet ensemble est vide dans le cas o\`u la valuation $\mu$ n'est pas bien sp\'ecifi\'ee, sinon c'est la r\'eunion d'un ensemble fini de boules ferm\'ees non vides de $\bar K$, qui correspondent aux diff\'erentes extensions de $\mu$ \`a $\bar K[x]$ ({\bf Th\'eor\`eme \ref{th:extensions}}). 

\vskip .2cm 

Dans la quatri\`eme partie nous montrons comment \`a partir d'une valuation $\bar\mu$ de $\bar K[x]$ nous pouvons construire la famille admise $\mathcal{A} = \bigl ( \mu _i \bigr ) _{i \in I}$ associ\'ee \`a la restriction $\mu$ de $\bar\mu$ \`a $K[x]$. 
Plus pr\'ecis\'ement nous construisons une famille d\'ecroissante de boules ferm\'ees $B_i$ de $\bar K$, telle que la valuation $\mu _i$ de la famille admise $\mathcal{A}$ soit la retriction \`a $K[x]$ de la valuation $\bar\mu _i$ de $\bar K[x]$ d\'efinie par la boule $B_i$. 

Alors que la construction de la famille admise $\mathcal{A} = \bigl ( \mu _i \bigr ) _{i \in I}$ se fait de \emph{mani\`ere croissante}, c'est-\`a-dire la valuation $\mu _i$ est construite \`a partir des valuations $\mu _j$ pour $j<i$, la construction des boules $B_i$, et par cons\'equent des valuations $\bar\mu _i$ se fait de \emph{mani\`ere d\'ecroissante}, c'est-\`a-dire la boule $B _i$ est construite \`a partir desboules $B _j$ pour $j>i$ ({\bf Th\'eor\`eme \ref{th:construction-famille}}). 

\vskip .2cm 
 
 Enfin dans l'annexe nous interpr\'etons les r\'esultats de Kaplansky sur les extensions imm\'ediates et les suites pseudo-convergentes \`a partir des propri\'et\'es des familles admises continues que nous avons d\'efinies pr\'ec\'edemment. 
  
\vskip .2cm

      \section{Rappels}    
%

Dans ce qui suit nous nous donnons une valuation $\nu$ sur un corps $K$ et toutes les valuations ou pseudo-valuations $\mu$ de l'anneau des polyn\^omes $K[x]$ que nous consid\'erons sont des prolongements de $\nu$. 
Nous nous donnons aussi un groupe totalement ordonn\'e $\tilde\Gamma$, contenant le groupe des ordres $\Gamma _{\nu}$ de la valuation $\nu$, et toutes les valuations ou pseudo-valuations $\mu$ de $K[x]$ ont leur groupe des ordres $\Gamma _{\mu}$ qui est un sous-groupe ordonn\'e de $\tilde\Gamma$.   

Pour toute valuation $\mu$ de $K[x]$ nous pouvons d\'efinir la notion de {\it polyn\^ome-cl\'e} $\phi$, et si $\phi$ est un polyn\^ome-cl\'e pour $\mu$ et si $\gamma$ est un \'el\'ement de $\tilde\Gamma$ v\'erifiant $\gamma > \mu (\phi )$, nous pouvons d\'efinir une nouvelle valuation $\mu '$ de $K[x]$, appel\'ee {\it valuation augment\'ee} associ\'ee au polyn\^ome-cl\'e $\phi$ et \`a la valeur $\gamma$ que nous notons $\mu ' = [ \mu \ ; \ \mu ' (\phi ) = \gamma ]$, de la mani\`ere suivante:  
   
pour tout polyn\^ome $f$ de $K[x]$, nous \'ecrivons le d\'eveloppement de $f$ selon les puissances de $\phi$, $f = g_m \phi ^m + \ldots + g_1 \phi + g_0$, o\`u les polyn\^omes $g_j$, $0 \leq j \leq m$, sont de degr\'e strictement inf\'erieur au degr\'e du polyn\^ome-cl\'e $\phi$, et nous avons:  
$$ \mu ' (f) \ = \ Inf \left ( \mu ( g_j ) + j \gamma \ ; \ 0 \leq j \leq m \right) \ .$$   

Nous pouvons d\'efinir aussi pour tout polyn\^ome unitaire $\phi$ de degr\'e un, $\phi = x-b$, et pour toute valeur $\gamma$ de $\tilde\Gamma$ une valuation $\mu$ de $K[x]$, que nous appelons encore {\it valuation augment\'ee} associ\'ee au polyn\^ome-cl\'e $\phi$ et \`a la valeur $\gamma$ que nous notons $\mu  = [ \nu \ ; \ \mu  (\phi ) = \gamma ]$, de la mani\`ere suivante:  

tout polyn\^ome $f$ de $K[x]$ s'\'ecrit de mani\`ere unique sous la forme $f= a_d \phi ^d + \ldots + a_1 \phi + a_0$, avec $a_j \in K$, et nous posons 
$$\mu (f) \ = \ Inf \left ( \nu (a_j) + j \gamma \ ; \ 0 \leq j \leq d \right ) \ .$$
Dans ce cas cette valuation $\mu$ est aussi not\'ee  $\omega _{(b,\gamma )}$ (cf. \cite{AP1}). 

\vskip .2cm 

Nous pouvons d\'efinir la notion de \emph{famille de valuations augment\'ees it\'er\'ees} comme une famille d\'enombrable $\bigl ( \mu _i \bigr ) _{i \in I}$ de valuations de $K[x]$, $I = \{ 1, \ldots ,n \}$ ou $I = \N ^*$, associ\'ee \`a une famille de polyn\^omes $\bigl (\phi _i\bigr )_{i \in I}$ et \`a une famille $\bigl (\gamma _i\bigr )_{i \in I}$ d'\'el\'ements de $\tilde\Gamma$, telle que chaque valuation $\mu _i$, $i>1$, est une valuation augment\'ee de la forme $\mu _i = [ \mu _{i-1} \ ; \ \mu _i (\phi _i) = \gamma _i ]$ et o\`u la famille des polyn\^omes-cl\'es $\bigl (\phi _i \bigr )$ v\'erifie les deux propri\'et\'es suivantes: pour tout $i>2$ nous avons deg$\phi _i \geq$ deg$\phi _{i-1}$ et les polyn\^omes $\phi _i$ et $\phi _{i-1}$ ne sont pas $\mu _{i-1}$-\'equivalents.    
Nous renvoyons aux articles \cite{McL 1}, \cite{McL 2}, et \cite {Va 1}, pour les d\'efinitions et les propri\'et\'es des polyn\^omes-cl\'es, des valuations augment\'ees et des familles de valuations augment\'ees it\'er\'ees.      

Nous d\'efinissons aussi la notion de \emph{famille admissible continue} comme une famille ${\mathcal C} = \bigl ( \mu _{\alpha} \bigr ) _{\alpha \in A}$ de valuations de $K[x]$, index\'ee par un ensemble totalement ordonn\'e $A$ sans plus grand \'el\'ement, associ\'ee \`a la famille de polyn\^omes-cl\'es $\bigl ( \phi _{\alpha} \bigr ) _{\alpha \in A}$ et \`a la famille de valeurs $\bigl ( \gamma _{\alpha} \bigr ) _{\alpha \in A}$. 
Par d\'efinition chaque valuation $\mu _{\alpha}$ est une valuation augment\'ee de la forme $\mu _{\alpha} = [ \mu  \ ; \ \mu _{\alpha} (\phi _{\alpha}) = \gamma _{\alpha} ]$, o\`u $\mu$ est une valuation de $K[x]$ donn\'ee, les polyn\^omes-cl\'es $\phi _{\alpha}$ sont tous de m\^eme degr\'e $d$ et les valeurs $\gamma _{\alpha}$ forment une famille croissante sans plus grand \'el\'ement dans $\tilde\Gamma$.   

Nous d\'efinissons l'ensemble 
$$\tilde\Phi \bigl ( \mathcal{C} \bigr ) = \tilde\Phi \left ( {\bigl (\mu _{\alpha}\bigr )}_{\alpha \in A} \right ) \ = \ \{ f \in K[x] \ | \ \mu _{\alpha} (f) < \mu _{\beta}(f) , \forall \alpha < \beta \in A \} \ ,$$  
nous supposons que cet ensemble est non vide, nous appelons $d _{ \mathcal{C}}$ le degr\'e minimal des polyn\^omes appartenant \`a $\tilde\Phi \bigl ( \mathcal{C} \bigr )$ et nous supposons aussi que nous avons l'in\'egalit\'e $d < d _{ \mathcal{C}}$, alors nous d\'efinissons l'ensemble  
$$\Phi \bigl ( \mathcal{C} \bigr ) = \Phi \left ( {\bigl (\mu _{\alpha} \bigr )}_{\alpha \in A} \right ) \ = \ \{ \phi \in \tilde\Phi \bigl ( \mathcal{C} \bigr ),\  \hbox{deg}\phi =  d _{ \mathcal{C}} \ \hbox{et $\phi$ unitaire} \  \} \ .$$  
Un polyn\^ome $\phi$ appartenant \`a $\Phi \bigl ( \mathcal{C} \bigr )$ est appel\'e un \emph{polyn\^ome-cl\'e-limite} pour la famille $\mathcal{C}$, et pour $\phi$ un polyn\^ome-cl\'e limite et $\gamma$ un \'el\'ement de $\tilde\Gamma$ v\'erifiant $\gamma > \mu _{\alpha}(\phi )$ pour tout $\alpha$ dans $A$, , nous pouvons d\'efinir une nouvelle valuation $\mu '$ de $K[x]$, appel\'ee \emph{valuation augment\'ee limite} pour $\mathcal{C}$ associ\'ee au polyn\^ome-cl\'e limite $\phi$ et \`a la valeur $\gamma$ que nous notons $\mu ' = \bigl [ ( \mu _{\alpha} ) _{\alpha \in A} \ ; \ \mu ' (\phi ) = \gamma \bigr ]$, de la mani\`ere suivante:  

pour tout polyn\^ome $f$ de $K[x]$, nous \'ecrivons le d\'eveloppement de $f$ selon les puissances de $\phi$, $f = g_m \phi ^m + \ldots + g_1 \phi + g_0$, o\`u les polyn\^omes $g_j$, $0 \leq j \leq m$, sont de degr\'e strictement inf\'erieur au degr\'e du polyn\^ome-cl\'e limite $\phi$, et nous posons:           
$$ \mu ' (f) \ = \ Inf \left (  \mu _A ( g_j ) + j \gamma \  ; \ 0 \leq j \leq m  \right ) \ ,$$    
o\`u nous posons $\mu _A(g) = Sup (\mu _{\alpha}(g) ; \alpha \in A)$ pour tout $g$ n'appartenant pas \`a $\tilde\Phi \bigl ( \mathcal{C} \bigr )$.  
Nous renvoyons \`a \cite{Va 1} pour les d\'efinitions pr\'ecises et les propri\'et\'es des poly\-n\^omes-cl\'es limites et des valuations augment\'ees limites.

\vskip .2cm 

\begin{remark}
Si nous prenons la valeur $\gamma = +\infty$, la valuation augment\'ee $\mu  = [ \nu \ ; \ \mu  (\phi ) = \gamma ]$ associ\'ee \`a un polyn\^ome-cl\'e $\phi$  ou la valuation augment\'ee limite $\mu ' = \bigl [ ( \mu _{\alpha} ) _{\alpha \in A} \ ; \ \mu ' (\phi ) = \gamma \bigr ]$ associ\'ee \`a un polyn\^ome-cl\'e limite $\phi$ est une pseudo-valuation de l'anneau $K[x]$ dont le noyau est l'id\'eal engendr\'e par le polyn\^ome $\phi$. 
\end{remark} 

\vskip .2cm 

\begin{definition}
Une \emph{famille admissible simple} $\mathcal S$ pour la valuation $\nu$ de $K$ est une famille de valuations $\bigl ( \mu _i \bigr ) _{i \in I}$ de $K[x]$ constitu\'ee d'une \emph{partie discr\`ete} ${\mathcal D}$ et d'une \emph{partie continue} ${\mathcal C}$,  

- la partie discr\`ete ${\mathcal D} = \bigl (\mu _l\bigr )_{l \in L}$ est une famille non vide de valuations augment\'ees it\'er\'ees de $K[x]$ telle que la famille de polyn\^omes-cl\'es $\bigl (\phi _l\bigr )_{l \in L}$ associ\'ee v\'erifie l'in\'egalit\'e stricte deg$\phi _l>$ deg$\phi _{l-1}$.

- la partie continue ${\mathcal C} = {(\mu _{\alpha})}_{\alpha \in A} $ est une famille admissible continue \'eventuellement vide; si elle est non vide la famille ${\mathcal D}$ est finie, le degr\'e $d$ des polyn\^omes-cl\'e $\phi _{\alpha}$ est \'egal au degr\'e du dernier polyn\^ome-cl\'e $\phi _n$ de la famille $\bigl (\phi _l \bigr ) _{l \in L}$ associ\'ee \`a ${\mathcal D}$, et pour tout $\alpha$ dans $A$, la valuation $\mu _{\alpha}$ est la valuation augment\'ee $\mu _{\alpha}= [\mu _n \ ; \ \mu _{\alpha} ( \phi _{\alpha} ) = \gamma _{\alpha} ]$. 
\end{definition} 

Si la partie discr\`ete ${\mathcal D}$ d'une famille admise simple ${\mathcal S}$ est constitu\'ee d'une seule valuation $\mu _1$, et si la partie continue ${\mathcal C}$ est non vide, nous pouvons toujours consid\'erer que la valuation $\mu _1$ appartient \`a la famille ${\mathcal C}$, nous \'ecrivons ${\mathcal S}={\mathcal C}$ et nous disons alors que la famille simple ${\mathcal S}$ est continue. 

\begin{definition}  
Une \emph{famille admissible} $\mathcal A$ pour la valuation $\nu$ de $K$ est une famille de valuations $\bigl ( \mu _i \bigr ) _{i \in I}$ de $K[x]$, obtenue comme r\'eunion de \emph{familles admissibles simples}     
$${\mathcal A} \ = \ \bigcup _{j \in J} {\mathcal S}^{(j)} \  \ = \  \bigcup _{j \in J}
\bigl ( {\mathcal D}^{(j)};{\mathcal C}^{(j)} \bigr ) \   \ , $$  
o\`u $J$ est un ensemble d\'enombrable, $J=\{ 1, \ldots , N\}$ ou $J = \N^*$, et nous d\'efinissons $J^*$ par $J^*=\{ 1, \ldots , N-1\}$ si $J$ est fini et par $J^*=J=\N^*$ sinon, v\'erifiant: 

- pour $j$ appartenant \`a $J^*$, la partie discr\`ete ${\mathcal D}^{(j)}  = \bigl (\mu _l^{(j)}\bigr )_{l \in L^{(j)}}$ est finie, la partie continue ${\mathcal C}^{(j)}  = {(\mu _{\alpha}^{(j)})}_{\alpha \in A^{(j)}} $ est non vide et la premi\`ere valuation $\mu _1^{(j+1)}$ de la famille simple ${\mathcal S}^{(j+1)}$ est une valuation augment\'ee limite pour la famille admissible continue ${\mathcal C}^{(j)}$; 

- la premi\`ere valuation $\mu ^{(1)}_1$ de la famille est la valuation associ\'ee \`a un polyn\^ome unitaire de degr\'e un, $\phi ^{(1)}_1 = x-a$, et \`a une valeur $\gamma ^{(1)}_1$, $\mu ^{(1)}_1 = [ \nu \ ; \ \mu ^{(1)}_1 ( \phi ^{(1)}_1) = \gamma ^{(1)}_1 ] = \omega _{(a, \gamma _1^{(1)})}$. 

Dans la suite, comme la valuation $\nu$ de $K$ est fix\'ee nous dirons simplement que ${\mathcal A}$ est une famille admissible de valuations de $K[x]$. 
\end{definition}

\vskip .2cm

Nous pouvons aussi \'ecrire la famille admissible ${\mathcal A}$ comme une famille index\'ee par un ensemble totalement ordonn\'e $I$,    
$${\mathcal A} \ = \ ( \mu _i )_{i \in I} \ , $$  
et l'ensemble $I$ peut \^etre d\'ecrit de la mani\`ere suivante:
pour tout $j$ dans $J$, nous munissons l'ensemble $B^{(j)} = L^{(j)} \sqcup A^{(j)}$ de l'ordre total induit par les ordres sur $L^{(j)}$ et sur $A^{(j)}$ et d\'efini par $l < \alpha$ pour tout $l \in L^{(j)}$ et tout $\alpha \in A^{(j)}$; et nous posons     
$$I \ = \ \bigl\{ (j,b) \ | \ j\in J \ \hbox{et} \ b \in B^{(j)} \bigr\}\ , $$   
muni de l'ordre lexicographique. 
L'ordre sur l'ensemble $I$ peut \^etre caract\'eris\'e par la relation suivante: $i<k$ dans $I$ si et seulement si pour tout polyn\^ome $f$ de $K[x]$ nous avons $\mu _i (f) \leq \mu _k(f)$ et il existe au moins un polyn\^ome $g$ avec $\mu _i(g)<\mu _k(g)$.       

La premi\`ere valuation $\mu _1$ de la famille $\cA$ est obtenue \`a partir de la valuation $\nu$ de $K$ gr\^ace \`a un polyn\^ome $\phi _1$ unitaire de degr\'e un et \`a une valeur $\gamma _1$.  
Nous consid\`ererons parfois que la valuation $\nu = \mu _0$ appartient \`a la famille $\cA$ et par abus de notation nous consid\`ererons que $0$ est le plus petit \'el\'ement de l'ensemble $I$.  
La valuation $\mu _1$ est ainsi consid\'er\'ee comme une valuation augment\'ee, d\'efinie par le polyn\^ome $\phi _1$.  

\vskip .2cm

A toute famille admissible ${\mathcal A}$ nous associons la famille des polyn\^omes-cl\'es ou polyn\^omes-cl\'es limites $\bigl ( \phi _i \bigr )_{i \in I}$, que nous appelons pour simplifier la famille des polyn\^omes-cl\'es, et la famille des valeurs $\bigl ( \gamma _i \bigr )_{i \in I}$. 

\vskip .2cm 

\begin{definition} 
 Une famille admissible $\mathcal{A} = ( \mu _i )_{i \in I}$ est une famille \emph{admise} si pour tout polyn\^ome $f$ dans $K[x]$ la famille $( \mu _i(f)) _{i \in I}$ admet un plus grand \'el\'ement dans le groupe $\Gamma$.
\end{definition}

\begin{definition}   
Une famille admissible ${\mathcal A} = \bigl ( \mu _i \bigr )_{i \in I}$ de valuations de $K[x]$ est dite \emph{compl\`ete} si l'ensemble $I$ poss\`ede un plus grand \'el\'ement $\bar\iota$, sinon la famille admissible ${\mathcal A}$ est dite \emph{ouverte}.   
\end{definition}   

\begin{remark}\label{rmq:famille_complete}
Une famille admissible ${\mathcal A}$ est compl\`ete uniquement dans le cas o\`u ${\mathcal A}$ est r\'eunion d'un nombre fini de familles simples et o\`u la derni\`ere famille simple ${\mathcal S}^{(N)}$ est discr\`ete finie,  ${\mathcal S}^{(N)} =\bigl (\mu _1^{(N)} ,\ldots ,\mu _{n_N}^{(N)} \bigr )$.      

Dans ce cas la famille ${\mathcal A}$ est admise et  la derni\`ere valuation $\mu _{\bar\iota} =  \mu _{n_N}^{(N)}$ de la famille peut \^etre une pseudo-valuation de $K[x]$. 
\end{remark}   

\begin{remark}\label{rmq:famille ouverte} 
Si la famille admissible $\mathcal{A} = ( \mu _i )_{i \in I}$ est ouverte, elle est admise si pour tout polyn\^ome $f$ il existe $i \in I$ tel que $\mu _i(f) = \mu _j(f)$ pour tout $j \geq i$.  
C'est le cas si la famille $\mathcal{A}$ est r\'eunion infinie de familles admissibles simples, ou si la famille $\mathcal{A}$ est r\'eunion de $N$ familles admissibles simples ${\mathcal A} = {\mathcal S}^{(1)} \cup \ldots \cup {\mathcal S}^{(N)}$,  telle que la derni\`ere famille simple ${\mathcal S}^{(N)}$ est une famille discr\`ete infinie, c'est-\`a-dire ${\mathcal S}^{(N)} = {\bigl ( \mu _l^{(N)}\bigr )}_{l \in L^{(N)}}$ avec $L^{(N)}$ infini, ou enfin si la famille simple ${\mathcal S}^{(N)}$ est de la forme ${\mathcal S}^{(N)} = \bigl ( {(\mu _l^{(N)})}_{l \in L^{(N)}}; {(\mu _{\alpha}^{(N)})}_{\alpha \in A^{(N)}}\bigr )$, avec $\tilde\Phi \bigl ( ( \mu _{\alpha}^{(N)}) _{\alpha \in A^{(N)}} \bigr ) = \emptyset$, c'est-\`a-dire telle que pour tout $f$ dans $K[x]$ il existe $\alpha < \beta$ dans $A^{(N)}$ avec $\mu _{\alpha} ^{(N)}(f) = \mu  _{\beta} ^{(N)}(f)$.  
\end{remark} 

\begin{definition} 
La famille admise ${\mathcal A} = ( \mu _i )_{i \in I}$ \emph{converge} vers la valuation ou pseudo-valuation $\mu$ de $K[x]$ d\'efinie pour tout polyn\^ome $f$ par 
$$\mu (f) \ = \ Sup \bigl ( \mu _i (f) \ ; \ i \in I \bigr ) \ .$$ 

Si l'ensemble $I$ admet un plus grand \'el\'ement $\bar \iota$, la limite de la famille est la valuation ou pseudo-valuation $\mu _{\bar \iota}$, sinon la limite est une valuation d\'efinie par $\mu (f) = \mu _i(f)$ pour $i$ assez grand dans $I$. 
\end{definition}

Nous avons une r\'eciproque au r\'esultat pr\'ec\'edent. 

\begin{theorem}\label{th1} (Th\'eor\`emes 2.4. et 2.5. de \cite{Va 1})
Soit $\mu$ une valuation ou pseudo-valuation de $K[x]$ prolongeant une valuation $\nu$ de $K$, alors il existe une famille admise de valuations de $K[x]$, not\'ee $\mathcal{A}(\mu )$ et appel\'ee \emph{famille admise associ\'ee \`a la valuation $\mu$} qui converge vers $\mu$. 
\end{theorem}

\begin{remark}\label{rmq:unicite}
La famille admise associ\'ee \`a une valuation $\mu$ n'est pas unique, mais est d\'etermin\'ee \`a \'equivalence pr\`es, o\`u deux familles admissibles ${\mathcal A} = \bigcup _{j \in J} {\mathcal S}^{(j)}$ et ${\mathcal A}' = \bigcup _{j \in J'} {\mathcal S}'^{(j)}$ sont dites \'equivalentes si $J=J'$, si les familles discr\`etes ${\mathcal D}^{(j)} = \bigl ( \mu ^{(j)}_i \bigr ) _{1 \leq i \leq n}$ et ${{\mathcal D}' }^{(j)} = \bigl ( {\mu '}^{(j)} _i \bigr ) _{1 \leq i \leq n'}$ co\"\i ncident jusqu'\`a l'avant-derni\`ere valuation, c'est-\`a-dire quand $n=n'$ et $\mu ^{(j)}_i = \mu '^{(j)}_i$ pour tout $i$, $1 \leq i \leq n-1$, et si les sous-familles continues ${\mathcal C}^{(j)}$ et ${\mathcal C}'^{(j)}$ co\"\i ncident asymptotiquement.  (cf. Proposition 2.9. de \cite{Va 2})
\end{remark}   

\vskip .2cm 

\begin{definition}   
Une valuation $\mu$ de $K[x]$ est dite \emph{bien sp\'ecifi\'ee} si la famille admise $\mathcal{A}(\mu )$ associ\'ee est compl\`ete.   
Dans ce cas la valuation $\mu$ est la derni\`ere valuation $\mu _{\bar\iota}$ de la famille $\mathcal{A}(\mu ) = ( \mu _i )_{i \in I}$.
\end{definition}   

Nous avons le r\'esultat suivant: 

\begin{proposition}  (Proposition 1.4 de \cite{Va 3})     
Les propositions suivantes sont \'equivalentes:                                        

1) La valuation $\mu$ est bien sp\'ecifi\'ee.              

2) La valuation $\mu$ n'est pas maximale pour la relation d'ordre $\leq$.       

3) La valuation $\mu$ admet un polyn\^ome-cl\'e.                      

4) La valuation $\mu$ peut \^etre obtenue comme valuation augment\'ee 
$$\mu = [ \mu _0 \ ; \ \mu (\phi )=\gamma ], $$             
ou comme valuation augment\'ee limite 
$$\mu = \bigl [ \bigl (\mu _{\alpha}\bigr ) _{\alpha \in A} \ ; \ \mu (\phi )=\gamma ] .$$
\end{proposition}

Si $\mu$ est une valuation bien sp\'ecifi\'ee, nous disons que le polyn\^ome $\phi _{\bar \iota}$ apparaissant comme dernier polyn\^ome de la famille $\bigl ( \phi _i \bigr ) _{i \in I}$ \emph{d\'efinit} la valuation ou pseudo-valuation $\mu$. 
Si $\mu$ est obtenue comme valuation augment\'e $\mu = [ \mu _0 \ ; \ \mu (\phi )=\gamma ]$, ou comme valuation augment\'ee limite, $\mu = \bigl [ \bigl (\mu _{\alpha}\bigr ) _{\alpha \in A} \ ; \ \mu (\phi )=\gamma ] $, nous pouvons en particulier choisir $\phi _{\bar\iota} = \phi$. 

Soit $\mu$ une valuation ou une pseudo-valuation de $K[x]$, alors toute valuation $\mu _i$ d'une famille admissible associ\'ee \`a $\mu$ est une valuation bien sp\'ecifi\'ee, d\'efinie par le polyn\^ome $\phi _i$. 

\vskip .2cm

Pour toute valuation ou pseudo-valuation $\mu$ de $K[x]$, les valuations $\mu _i$ appartenant \`a une famille admise associ\'ee \`a $\mu$ sont d\'efinies de mani\`ere essentiellement unique (cf. remarque \ref{rmq:unicite}), en particulier quand $\mu$ est bien sp\'ecifi\'ee, si $\mu$ est une valuation augment\'ee, $\mu = [ \mu _0 \ ; \ \mu (\phi )=\gamma ]$, la valuation $\mu _0$ est d\'efinie de mani\`ere unique, et si $\mu$ est une valuation augment\'ee limite, $\mu = \bigl [ \bigl (\mu _{\alpha}\bigr ) _{\alpha \in A} \ ; \ \mu (\phi )=\gamma ]$, la famille $\bigl (\mu _{\alpha}\bigr ) _{\alpha \in A}$ est bien d\'efinie \emph{asymptotiquement}. 

Dans la suite nous noterons alors 
$$\mu = [ \mu _{\sharp} \ ; \ \mu (\phi )=\gamma ], $$  
o\`u $\mu _{\sharp}$ est la valuation $\mu _0$, resp. une famille continue de valuations $\bigl (\mu _{\alpha}\bigr ) _{\alpha \in A}$, et o\`u $\phi$ est un polyn\^ome-cl\'e, resp. un polyn\^ome-cl\'e limite, pour $\mu _{\sharp}$.
En g\'en\'eral le polyn\^ome $\phi$ n'est pas d\'efini de mani\`ere unique, en fait les polyn\^omes $\phi$ et $\psi$ d\'efinissent la m\^eme valuation $\mu$ si et seulement si ce sont des polyn\^omes unitaires de m\^eme degr\'e v\'erifiant $\mu _{\sharp}(\phi - \psi) \geq \gamma$, o\`u nous posons $\mu _{\sharp}(f) = \mu _A(f) = Sup (\mu _{\alpha}(f) ; \alpha \in A)$ dans le cas o\`u $\mu _{\sharp}$ est la famille $\bigl (\mu _{\alpha}\bigr ) _{\alpha \in A}$ (\cite{Va 2}). 

\vskip .2cm

\begin{remark}
Comme un polyn\^ome-cl\'e est irr\'eductible, dans le cas o\`u le corps $\bar K$ est alg\'ebri\- quement clos, une valuation bien sp\'ecifi\'ee $\mu$ est de la forme 
$$\mu = [ \nu \ ; \ \mu (\phi ) = \delta ] , $$ 
avec $\phi(x) = x-a$. 
Cela correspond \`a la valuation associ\'ee \`a la paire minimale $(a,\delta )$ not\'ee $\omega _{(a,\delta )}$ d\'efinie dans \cite{AP1}.  
\end{remark}

\vskip .2cm

      \section{Groupe des valeurs et alg\`ebre gradu\'ee}    
%

Soit $\mu$ une valuation sur un corps $K$ de groupe des valeurs $\Gamma _{\mu}$, pour tout sous-anneau $A$ de $K$ et pour tout $\gamma$ dans $\bar \Gamma _{\mu} = \Gamma _{\mu} \cup \{ + \infty \}$, nous d\'efinissons les groupes ${\cP}_{\gamma}(A) = \{ x \in A \ | \ \mu (x) \geq \gamma \}$ et ${\cP}_{\gamma}^+ (A) = \{ x \in A\ | \ \mu (f) >\gamma \}$, et l'alg\`ebre gradu\'ee ${\rm gr}_{\mu}A$ associ\'ee \`a la valuation $\mu$ par : 
$${\rm gr}_{\mu} A \ = \ \bigoplus _{\gamma \in \bar\Gamma} {\cP}_{\gamma} (A)/ {\cP}_{\gamma}^+ (A)\ .$$
Nous notons $H_{\mu}$ l'application de $A$ dans ${\rm gr}_{\mu}A$ qui \`a tout \'el\'ement $x$ de $A$ avec $\mu (x) = \gamma$ associe l'image de $x$ dans ${\cP}_{\gamma} (A)/ {\cP}_{\gamma}^+ (A)$, et nous notons $\Delta _{\mu}(A)$ la partie homog\`ene de degr\'e $0$, $\Delta _{\mu}(A) = {\cP}_{0} (A)/ {\cP}_{0}^+ (A)$.  

En particulier pour $A=K$ les groupes ${\cP}_{0} (K)$ et ${\cP}_{0}^+ (K)$ sont respectivement l'anneau $V_{\mu}$ de la valuation et son id\'eal maximal, $\Delta _{\mu}(K)$ est \'egal \`a son corps r\'esiduel $\kappa _{\mu}$ et l'alg\`ebre gradu\'ee ${\rm gr}_{\mu}K$ est \emph{simple}, i.e. tout \'el\'ement homog\`ene non nul admet un inverse.  
Plus g\'en\'eralement si $K$ est le corps des fractions de $A$, l'alg\`ebre gradu\'ee  ${\rm gr}_{\mu}K$ est l'alg\`ebre gradu\'ee simple engendr\'ee par ${\rm gr}_{\mu}A$ et le corps r\'esiduel $\kappa _{\mu}$ est le corps des fractions de l'anneau $\Delta _{\mu} (A)$. 

\vskip .2cm

Soit $\cA = \bigl ( \mu _i \bigr ) _{i \in I}$ une famille admissible de valuations, et pour tout $i \in I$ nous appelons $\Gamma _{\mu _i}$ le groupe des ordres de la valuation $\mu _i$. 

Si $\mu _k$ et $\mu _l$ sont deux valuations appartenant \`a la m\^eme sous-famille simple $\cS$ de la famille $\cA$ telles que $\mu _l$ est obtenue comme valuation augment\'ee $\mu _l=[\mu _k \ ; \ \mu _l(\phi _l)=\gamma _l]$, nous disons que $(\mu _k,\mu _l)$ forment {\it un couple de valuations successives} de la famille.
Le groupe des ordres $\Gamma _{\mu _l}$ de la valuation $\mu _l$ est \'egal \`a $\Gamma _{\mu _l} = \Gamma _{\mu _k} \oplus \Z \gamma _l$, d'o\`u l'\'egalit\'e 
$$ [\Gamma _l : \Gamma _k] \ = \ \tau _l \,$$  
o\`u $\tau _l$ est le plus petit entier $t>0$ tel que $t\gamma _l$ appartienne \`a $\Gamma _{\mu _k}$ si $\gamma _l$ appartient \`a $\Gamma _{\mu _l} \otimes _{\Z}\Q$, et o\`u $\tau _l$ est $+\infty$ sinon.   
Remarquons que la valuation $\mu _l$ admet un polyn\^ome-cl\'e qui n'est pas $\mu _l$-\'equivalent au polyn\^ome $\phi _l$ si et seulement si la valeur $\gamma _l$ appartient au groupe $\Gamma _k \otimes _{\Z} \Q$, en particulier si $\gamma _l$ n'appartient pas \`a $\Gamma _{\mu _l} \otimes _{\Z}\Q$ la valuation $\mu _l$ est la derni\`ere valuation de la famille admissible $\cA$. 

Comme les valuations $\mu _k$ et $\mu _l$ v\'erifient $\mu _k (f) \leq \mu _l (f)$ pour tout $f$ dans $K[x]$ nous avons une application naturelle $g \colon {\rm gr}_{\mu _k}K[x] \to {\rm gr}_{\mu _l}K[x]$, et celle-ci induit un isomorphisme    
$$G \colon \bigl ( {\rm gr} _{\mu _k} K[x] / (H _{\mu _k} ( \phi _l)) \bigr ) [T] \ 
\llrightarrow \ {\rm gr} _{\mu _l} K[x] \ , $$   
qui envoie $T$ sur $G(T) = H _{\mu _l} ( \phi _l)$ (cf. \cite{Va 1}).

Rappelons qu'il existe $q_k$ et $q'_k$ dans $K[x]$ v\'erifiant $q_k q'_k$ $\mu _k$-\'equivalent \`a $1$ et $\mu _k (q _k) = - \mu _k (q' _k) = \mu _k ( \phi _l)$, et nous posons $\varphi _l = H_{\mu _k} (q' _k \phi _l)$. 
De plus si $\gamma _l$ appartient \`a $\Gamma _{\mu _k} \otimes _{\Z} {\Q}$, il existe $p _l$ et $p' _l$ v\'erifiant $p _l p' _l$ $\mu _l$-\'equivalent \`a $1$ et $\mu _l (p _l) = - \mu _l(p' _l) =  \tau _l\gamma _l$ (cf. \cite{Va 3}).
Alors le noyau de la composante de degr\'e 0 de l'application $g$, $g_0 \colon \Delta _{\mu _k} \to \Delta _{\mu _l}$, est l'id\'eal engendr\'e par $\varphi _l$, et nous avons:    

\noindent - si $\gamma _l$ n'appartient pas \`a $\Gamma _{\mu _k} \otimes _{\Z} {\Q}$   
$$\bigl (\Delta _{\mu _k}/ (\varphi _l)\bigr ) \ \buildrel\sim\over\llrightarrow \ 
\Delta _{\mu _l} \ ,$$     

\noindent - si $\gamma _l$ appartient \`a $\Gamma _{\mu _k} \otimes _{\Z} {\Q}$ 
$$\bigl (\Delta _{\mu _k}/ (\varphi _l)\bigr ) [S _l] \ \buildrel\sim\over\llrightarrow \ 
\Delta _{\mu _l} \ ,$$   
avec $S _l=H _{\mu _l}\bigl ( p' _l {\phi _l}^{\tau _l}\bigr )$ 
(cf. \cite{Va 1} Remarque 1.5).    

\vskip 0.2cm

Si $\mu _l$ est la valuation augment\'ee limite d'une famille continue ${\cC}=\bigl ( \mu _{\alpha}\bigr )_{\alpha \in A}$, associ\'ee au polyn\^ome cl\'e-limite $\phi _l$, $\mu _l = \bigl [ \bigl (\mu _{\alpha}\bigr ) _{\alpha \in A} \ ;\ \mu _l(\phi _l)=\gamma _l\bigr ]$, nous d\'efinissons l'alg\`ebre gradu\'ee 
${\bf gr _A} \ = \ {\rm gr}_{\mu _{\alpha}}K[x]/(H_{\mu _{\alpha}} (\phi _{\beta}))$ 
qui ne d\'epend pas du couple $\alpha < \beta$ dans $A$, et l'application naturelle de ${\rm gr}_{\mu _{\alpha}}K[x]$ dans ${\rm gr}_{\mu _l}K[x]$ induit un isomorphisme d'alg\`ebres gradu\'ees:
$$Q \colon {\bf gr _A} [T] \ \buildrel\sim\over\llrightarrow \ {\rm gr}_{\mu _l}K[x] \ ,$$
qui envoie $T$ sur $Q(T) = H_{\mu _l}(\phi _l)$.
Nous appelons ${\bf \Delta _A}$ la composante de degr\'e $0$ de ${\bf gr _A}$, cet anneau est isomorphe \`a $\Delta _{\mu _{\beta}}/(\varphi _{\alpha})$ 
o\`u $(\mu _{\alpha},\mu _{\beta})$ est un couple de valuations successives de $\cA$ appartenant \`a ${\cC}$, avec 
$\varphi _{\alpha} = H_{\mu _{\alpha}} \bigl ( q' _{\alpha} \phi _{\beta} \bigr )$.

Tous les groupes de valuation $\Gamma _{\mu _{\alpha}}$ sont \'egaux et nous notons ce groupe ${\bf \Gamma _A}$. 
Pour tout $\gamma$ dans ${\bf \Gamma _A}$ il existe $p$ et $p'= p'(\gamma)$ dans $K[x]$ v\'erifiant $pp'(\gamma) \equiv{\mu _{\alpha}}1$ et $\mu _{\alpha}(p) = - \mu _{\alpha} (p'(\gamma ))=\gamma$, pour $\alpha \in A$. 

De plus si $\gamma _l$ appartient \`a ${\bf \Gamma _A} \otimes _{\Z} \Q$ et si nous appelons comme pr\'ec\'edemment $\tau _l$ le plus petit entier $t >0$ tel que $t\gamma _l$ appartienne \`a ${\bf \Gamma _A}$ il existe $p$ et $p'$ dans $K[x]$ tels que $pp'$ soit $\mu _{\alpha}$-\'equivalent \`a $1$ pour $\alpha$ suffisamment grand et tels que $\mu _{\alpha} (p') = -\tau _l\gamma _l$ (cf. \cite{Va 3} Proposition 2.2).    

Alors le morphisme $Q$ induit un isomorphisme en degr\'e $0$:  

\noindent - si $\gamma _l$ n'appartient pas \`a ${\bf \Gamma _A} \otimes _{\Z} \Q$     

$$Q_0 \colon {\bf \Delta _A} \ \buildrel\sim\over\llrightarrow \ \Delta _{\mu _l} \ ,$$

\noindent - si $\gamma _l$ appartient \`a ${\bf \Gamma _A} \otimes _{\Z} \Q$    

$$Q_0 \colon {\bf \Delta _A}[S] \ \buildrel\sim\over\llrightarrow \ \Delta _{\mu _l} \ ,$$
qui envoie $S$ sur $H_{\mu _l}(p' {\phi _l}^{\tau _l})$.    
  
\vskip 0.2cm

\begin{remark} 
Soit $\mu _l$ une valuation de la famille $\cA$, nous notons $\Gamma _{\sharp}$ le groupe des ordres $\Gamma _{\mu _k}$ de la valuation $\mu _k$ si $\mu _l$ est obtenue comme valuation augment\'ee,  
$\mu _l = [ \mu _k \ ; \ \mu _l(\phi _l) = \gamma _l ]$, 
ou le groupe des ordres $\bf \Gamma _A$ si la valuation $\mu _l$ est obtenue comme valuation augment\'ee limite, 
$\mu _l = \bigr [ \bigr ( \mu _{\alpha} \bigl ) _{\alpha \in A} \ ; \ \mu _l(\phi _l) = \gamma _l \bigr ]$.  

Si $\mu _l$ n'est pas la derni\`ere valuation de la famille $\cA$, 
il existe une valuation $\mu _m$ telle que $(\mu _l,\mu _m)$ est un couple de valuations 
successives, et nous \'ecrivons le polyn\^ome-cl\'e $\phi _m$ sous la forme
$\phi _m = {\phi _l}^{r_l} + \ldots + g_0$, nous avons  $r_l \gamma _l \in \Gamma _{\sharp}$, 
en particulier $\gamma _l$ appartient \`a $\Gamma _{\sharp} \otimes _{\Z} \Q$
et nous pouvons d\'efinir l'entier $s_l$ par $r_l = \tau _l s_l$.  
\end{remark}

\begin{proposition}\label{familledecorps}  
Il existe une famille croissante de corps $\bigl ( F_k \bigr ) _{k\in I^*}$, avec 
$F_0$ \'egal au corps r\'esiduel $\kappa _{\nu}$ de la valuation $\nu$ de $K$, 
telle que pour tout couple $(\mu _k,\mu _l)$ de valuations successives de $\cA$ nous avons: 

\noindent - si $\gamma _l$ appartient \`a $\Gamma _{\mu _k} \otimes _{\Z}\Q$ 
$$\Delta _{\mu _l} \ = \ F_k[S_l] \ , \quad   
\hbox{avec} \ S_l = H_{\mu _l} \bigl ( p'_l \phi _l^{\tau _l} \bigr ) \ ;$$   

\noindent - si $\gamma _l$ n'appartient pas \`a $\Gamma _{\mu _k} \otimes _{\Z}\Q$ 
$$\Delta _{\mu _l} \ = \ F_k \ . $$ 

De plus si $l$ appartient \`a $I^*$, $F_l$ est un extension finie de $F_k$ de degr\'e $s_l$, et pour $l$ tel que la valuation $\mu _l$ appartienne \`a une famille continue $\cC = \bigl ( \mu _{\alpha} \bigr ) _{\alpha \in A}$, le corps $F_l$ est isomorphe \`a ${\bf \Delta _A}$. 

En particulier tous les corps $F_l$ sont des extensions alg\'ebriques du corps r\'esiduel $\kappa _{\nu}$, et si la famille $\cA$ est constitu\'ee d'un nombre fini de sous-familles simples, tous les corps $F_l$ sont des extensions finies de $\kappa _{\nu}$. 
\end{proposition}    

\begin{preuve}   
La proposition est une g\'en\'eralisation du r\'esultat de MacLane (cf. \cite{McL 1} 
Theorem 12.1 et \cite{Va 1} Th\'eor\`eme 1.12) et se d\'emontre par r\'ecurrence (cf. \cite{Va 3}) . 

\hfill$\Box$ 
\end{preuve}    

\begin{remark}  
Nous avons montr\'e de plus que si $\mu _k$ est la premi\`ere valuation $\mu _1^{(j)}$ 
d'une sous-famille simple ${\cS}^{(j)}$, et si nous notons $F_0^{(j)}$ le corps tel que 
$\Delta _{\mu _k}$ soit \'egal \`a $F_0^{(j)}[S]$, alors $F_k$ est une extension 
alg\'ebrique finie de $F_0^{(j)}$ de degr\'e $s_k$.  
En effet nous avons $S=H_{\mu _k} \bigl ( p'_k{\phi _k}^{\tau _k}\bigr )$ et le corps 
$F_k$ est \'egal \`a $\Delta _{\mu _k} / (\varphi _l)$ o\`u 
$\varphi _l=H_{\mu _k} (q'_l\phi _l)$ avec $\phi _l = {\phi _k}^{\tau _ks_k}+\ldots + g_0$.  
\end{remark}   

\vskip .2cm

\begin{proposition} (Proposition 2.3 de \cite{Va 3})\label{prop:algebresimple} 
Soit $\mu$ une valuation de l'anneau des polyn\^omes $K[x]$, alors l'alg\`ebre gradu\'ee associ\'ee ${\rm gr}_{\mu}K[x]$ est de la forme suivante: 

 i) si la valuation $\mu$ n'est pas bien sp\'ecifi\'ee  
$${\rm gr}_{\mu}K[x] \ = \ {G^{(0)}} \ ,$$
o\`u ${G^{(0)}}$ est une alg\`ebre gradu\'ee simple, c'est-\`a-dire telle que tout \'el\'ement homog\`ene non nul admette un inverse; 

ii) si la valuation $\mu$ est bien sp\'ecifi\'ee  
$${\rm gr}_{\mu}K[x] \ = \ {G^{(0)}} [T]\ ,$$
o\`u ${G^{(0)}}$ est une alg\`ebre gradu\'ee simple et $T$ est l'image $H_{\mu}(\phi )$ du polyn\^ome $\phi$ d\'efinissant la valuation $\mu$. 

De plus un \'el\'ement homog\`ene $\psi$ de ${\rm gr}_{\mu}K[x]$ est irr\'eductible si et seulement si il existe $f$ polyn\^ome-cl\'e pour la valuation $\mu$ dans $K[x]$ et $\varepsilon$ \'el\'ement homog\`ene inversible de ${\rm gr}_{\mu}K[x]$ tels que $\varepsilon \psi$ soit \'egal \`a l'image $H_{\mu}(f)$ de $f$ dans ${\rm gr}_{\mu}K[x]$.  
\end{proposition}

\begin{proposition}\label{prop:egalite-d-Abhyankar}
La valuation $\mu$ de $K[x]$ est bien sp\'ecifi\'ee si et seulement si l’extension $(K(x), \mu )/(K, \nu)$ de corps valu\'es v\'erifie l'\'egalit\'e  d’Abhyankar:
$$dim.alg._K K(x) = dim.alg._{\kappa _{\nu}} \kappa _{\mu} + rang.rat. \Gamma _{\mu} / \Gamma _{\nu}= 1.$$
\end{proposition} 

\begin{preuve} 
Rappelons que le corps r\'esiduel $\kappa _{\mu}$ est \'egal au corps des fractions de la partie homog\`ene de degr\'e $0$, $\Delta _{\mu}$ de l'alg\`ebre gradu\'ee ${\rm gr}_{\mu}K[x]$.  

Dans le cas o\`u la valuation $\mu$ est obtenue comme valuation augment\'ee, $\mu = [ \mu _{\sharp} \ ; \ \mu (\phi ) = \gamma ]$, l'alg\`ebre gradu\'ee  ${\rm gr}_{\mu}K[x]$ est de la forme $G^{(0)} [T]$ o\`u l'alg\`ebre simple $G^{(0)}$ est isomorphe \`a l'alg\`ebre quotient  ${\rm gr}_{\mu _{\sharp}}K[x] /( H_{\mu _{\sharp}}( \phi ))$ et o\`u $T = H_{\mu }( \phi )$. 

Si $\gamma$ n'appartient pas \`a $\Gamma _{\mu _{\sharp}} \otimes _{\Z} \Q$, la partie homog\`ene de degr\'e $0$, $\Delta _{\mu }$, est isomorphe \`a $\Delta _{\mu _{\sharp}}  / (\varphi _{\sharp})$, c'est-\`a-dire au corps $F _{\sharp}$, nous en d\'eduisons que le corps r\'esiduel $\kappa _{\mu }$ de la valuation $\mu$ est isomorphe \`a  $F _{\sharp}$, par cons\'equent est une extension alg\'ebrique finie du corps r\'esiduel $\kappa _{\nu }$. 

Si $\gamma$ appartient \`a $\Gamma _{\mu _{\sharp}} \otimes _{\Z} \Q$, la partie homog\`ene de degr\'e $0$, $\Delta _{\mu }$, est isomorphe \`a $\Delta _{\mu _{\sharp}}  / (\varphi _{\sharp}) [S]$, avec $S = H_{\mu }( p' \phi ^{\tau } )$ o\`u $\tau$ est \'egal \`a $[\Gamma _{\mu} : \Gamma _{\nu}]$, et le corps r\'esiduel $\kappa _{\mu }$ de la valuation $\mu$ est isomorphe \`a  $F _{\sharp}(S)$, par cons\'equent est une extension transcendante de degr\'e $1$ du corps r\'esiduel $\kappa _{\nu }$. 

Dans le cas o\`u la valuation $\mu$ est obtenue comme valuation augment\'ee limite, $\mu = \bigl [ \bigl (\mu _{\alpha}  \bigr ) \ ; \ \mu (\phi ) = \gamma \bigr ]$, nous avons un r\'esultat analogue. 

Si $\gamma$ n'appartient pas \`a ${\bf \Gamma _A} \otimes _{\Z} \Q$, la partie homog\`ene $\Delta _{\mu }$ est isomorphe \`a ${\bf \Delta _A}$, le corps r\'esiduel $\kappa _{\mu }$ de la valuation $\mu$ est isomorphe au corps ${\bf \Delta _A}$ et est donc une extension alg\'ebrique finie du corps r\'esiduel $\kappa _{\nu }$. 

Si $\gamma$ appartient \`a ${\bf \Gamma _A} \otimes _{\Z} \Q$, la partie homog\`ene $\Delta _{\mu }$ est isomorphe \`a ${\bf \Delta _A} [S]$, avec $S = H_{\mu _{\alpha}}( p' \phi ^{\tau } )$ o\`u $\tau$ est \'egal \`a $[\Gamma _{\mu} : {\bf \Gamma _A}]$, et le corps r\'esiduel $\kappa _{\mu }$ de la valuation $\mu$ est isomorphe \`a  ${\bf \Delta _A}(S)$, par cons\'equent est une extension transcendante de degr\'e $1$ du corps r\'esiduel $\kappa _{\nu }$. 

\vskip .2cm 

Si la valuation $\mu$ n'est pas bien sp\'ecifi\'ee, chacune des valuations $\mu _l$ de la famille admise $\cA$ associ\'ee \`a la valuation $\mu$ a un groupe des ordres $\Gamma _{\mu _l}$ qui est une extension finie du groupe $\Gamma _{\nu}$, donc le groupe $\Gamma _{\mu}$, r\'eunion des groupes  $\Gamma _{\mu _l}$ a m\^eme rang rationnel que le groupe $\Gamma _{\nu}$. 

Le corps r\'esiduel $\kappa _{\mu}$ est la r\'eunion des corps $F_l$, extensions finies de $\kappa _{\nu}$, donc une extension alg\'ebrique du corps r\'esiduel $\kappa _{\nu}$. 
\hfill$\Box$ 
\end{preuve}

Nous pouvons d\'eduire de ce qui pr\'ec\`ede le r\'esultat suivant, qui r\'epond \`a une question pos\'ee par Nagata (cf. \cite{Na}) et a \'et\'e r\'esolue par J. Ohm (\cite{Oh}). 
Rappelons que nous disons qu'une extension de corps $l/k$ est \emph{r\'egl\'ee} s'il existe $k \subset k_1 \subset l$ avec $l/k_1$ extension transcendante pure de degr\'e $1$ et $k_1/k$ extension alg\'ebrique finie. 

\begin{corollary}  (\emph{\bf The ruled residue conjecture})
Soit $(K(x), \mu )/(K, \nu)$ une extension de corps valu\'es, alors le corps r\'esiduel $\kappa _{\mu}$ est une extension alg\'ebrique ou r\'egl\'ee du corps r\'esiduel $\kappa _{\nu}$. 
\end{corollary}  

\vskip .2cm

Le rang $rg(\mu )$ de la valuation $\mu$ est compris entre $rg(\nu )$ et $rg(\nu ) +1$, la valuation $\mu$ a le m\^eme rang que la valuation $\nu$ si $\gamma$ appartient au groupe $\Gamma _{\nu} \otimes _{\Z} \R$, sinon la valeur $\gamma$ appartient \`a un groupe totalement ordonn\'e $\tilde\Gamma$ qui contient $\Gamma _{\nu}$ comme sous groupe isol\'e et $\gamma$ v\'erifie $\gamma > \delta$ pour tout $\delta$ dans $\Gamma _{\nu}$. 

Dans ce ce dernier cas la valuation $\mu$ est essentiellement unique, c'est-\`a-dire que si nous nous donnons un polyn\^ome $\phi$ qui est polyn\^ome-cl\'e pour une valuation $\mu _{\sharp}$ ou polyn\^ome-cl\'e limite pour une famille  de valuation $\bigl ( \mu _{\alpha} \bigr ) _{\alpha \in A}$, la valuation bien sp\'ecifi\'ee $\mu$ d\'efinie par le polyn\^ome $\phi$ et la valeur $\gamma$ est ind\'ependante \`a \'equivalence pr\`es de la valeur $\gamma$ choisie dans $\tilde\Gamma \setminus \Gamma _{\nu}$.  

De plus le polyn\^ome $\phi$ qui d\'efinit une valuation $\mu$ de rang $rg(\mu )=rg(\nu )+1$ est unique. En effet si deux polyn\^omes $\phi$ et $\psi$ d\'efinissent la m\^eme valuation $\mu$ comme valuation augment\'ee ou comme valuation augment\'ee limite avec la valeur $\gamma$, nous avons l'in\'egalit\'e $\mu (\psi - \phi) \geq  \gamma$. 

\begin{remark}\label{rmq:valuation=pseudo-valuation}
Si nous prenons $\gamma =+\infty$ nous trouvons une pseudo-valuation de $K[x]$ dont le noyau est \'egal \`a l'id\'eal engendr\'e par $\phi$. 
Il y a une bijection entre l'ensemble des valuations $\mu$ de $K[x]$ de rang $rg(\mu ) = rg(\nu )+1$ et l'ensemble des pseudo-valuations de $K[x]$ de noyau non trivial, et l'\'etude des valuations de rang $rg(\mu ) = rg(\nu )+1$ d\'efinies par le polyn\^ome $\phi$ est \'equivalente \`a l'\'etude des pseudo-valuations de noyau $(\phi )$, c'est-\`a-dire \`a l'\'etude des valuations de l'extension $L=K[x]/(\phi )$ de $K$ qui prolongent $\nu$. 
\end{remark}
\vskip .2cm

\begin{proposition}\label{prop:generateurs} 
Soit $\mu$ une valuation bien sp\'ecifi\'ee de $K[x]$ d\'efinie par le polyn\^ome $\phi$,  et soit ${\rm gr}_{\mu}K[x]=G^{(0)} [T]$ l'alg\`ebre gradu\'ee associ\'ee avec $T = H_{\mu}(\phi)$, alors si $\psi$ est un autre polyn\^ome qui d\'efinit la valuation $\mu$ nous avons $S=H_{\mu}(\psi )$ qui est \'egal \`a $T$ ou \`a $T -h$ avec $h \in G^{(0)}$ de valuation $\mu(h ) = \mu (\phi ) = \mu(\psi )$. 

R\'eciproquement tout g\'en\'erateur homog\`ene $S$ de l'alg\`ebre gradu\'ee ${\rm gr}_{\mu}K[x]$ sur l'alg\`ebre simple ${G^{(0)}}$ est de la forme $S=T-h$ avec $h \in G^{(0)}$ de degr\'e $\mu(h ) = \mu (\phi ) = \mu(\psi )$, et il existe un polyn\^ome $\psi$ dans $K[x]$ qui d\'efinit la valuation $\mu$ avec $H_{\mu}(\psi ) =S$. 
\end{proposition} 

\begin{preuve} 
Si la valuation $\mu$ est obtenue comme valuation augment\'ee $\mu = [ \mu ' \ ; \ \mu (\phi ) = \gamma ]$ c'est une cons\'equence du r\'esultat suivant: 

deux valuations augment\'ees $\mu _1$ et $\mu _2$ d'une m\^eme valuation $\mu$ d\'efinies respectivement par des polyn\^omes-cl\'es $\phi _1$ et $\phi _2$ et des valeurs $\gamma _1$ et $\gamma _2$ sont \'egales si et seulement si $\gamma _1 = \gamma _2$ et si les polyn\^omes $\phi _1$ et $\phi _2$ ont m\^eme degr\'e et v\'erifient $\mu (\phi _1 - \phi _2) \geq \gamma _1$ (Proposition 1.2. de \cite{Va 2}). 

Si la valuation $\mu$ est obtenue comme valuation augment\'ee limite $\mu = \bigl [ \bigl (\mu _{\alpha}\bigr ) _{\alpha \in A} \ ; \ \mu (\phi )=\gamma ]$ c'est une cons\'equence du r\'esultat analogue: 

deux valuations augment\'ees limites  $\mu _1$ et $\mu _2$ d'une m\^eme  famille admissible continue $\mathcal{C} =  \bigl ( \mu _{\alpha} \bigr ) _{\alpha \in A}$ d\'efinies respectivement par des polyn\^omes-cl\'es limites $\phi _1$ et $\phi _2$ et des valeurs $\gamma _1$ et $\gamma _2$ sont \'egales si et seulement si $\gamma _1 = \gamma _2$ et si les polyn\^omes $\phi _1$ et $\phi _2$ ont m\^eme degr\'e et v\'erifient $\mu _A (\phi _1 - \phi _2) \geq \gamma _1$ (Proposition 1.4. de \cite{Va 2}). 
\hfill \tf    
\end{preuve}  

\vskip .2cm 

Nous rappelons aussi le r\'esultat suivant, qui est une cons\'equence de la proposition pr\'ec\'edente, mais qui peut se d\'emontrer aussi directement \`a   partir des propositions 1.2 et 1.4 de \cite{Va 2}. 

\begin{proposition}\label{prop:polynome-definissant-la-valuation} 
Soit $\mu$ une valuation bien sp\'ecifi\'ee de $K[x]$ d\'efinie par le polyn\^ome $\phi$, et soit $\psi$ un polyn\^ome unitaire de $K[x]$ v\'erifiant ${\rm deg}\ \psi = {\rm deg}\ \phi$ et $\mu (\psi ) = \mu (\phi )$, alors le polyn\^ome $\psi$ d\'efinit la valuation $\mu$. 
\end{proposition}
\hfill \tf  

\vskip .2cm 

Dans la suite de ce paragraphe nous allons supposer que le corps $K$ est alg\'ebriquement clos, alors pour toute valuation $\nu$ de $K$ le corps r\'esiduel $\kappa _{\nu}$ est aussi alg\'ebriquement clos et le groupe des valeurs $\Gamma _{\nu}$ est divisible. 
De plus dans le cas o\`u $K$ est alg\'ebriquement clos, les \'el\'ements irr\'eductibles de l'anneau $K[x]$ sont les polyn\^omes de degr\'e $1$, nous en d\'eduisons que toute valuation $\mu$ de $K[x]$ est d\'efinie enti\`erement par les valeurs $\mu (x-b)$, pour $b\in K$, et que tout polyn\^ome $f$ de degr\'e plus grand que $2$ ne peut pas \^etre un polyn\^ome-cl\'e ou un polyn\^ome-cl\'e limite. 

Nous rappelons que pour tout corps $L$, pour trouver la famille admise $\cA$ associ\'ee \`a une valuation $\mu$ de $L[x]$ le premier pas est de consid\'erer l'ensemble $\Lambda _{\mu} = \left\{ \mu (x-b) \ | \ b \in L \right\}$. 
Si cet ensemble a un plus grand \'el\'ement $\delta$ nous choisissons un polyn\^ome $\phi = x-a$  pour lequel cette valeur est atteinte et la premi\`ere valuation $ \mu _1$ de la famille est la valuation associ\'ee, c'est-\`a-dire la valuation $\mu _1= \omega _{(a,\delta )}$. Alors soit la valuation $\mu _1$ est la valuation $\mu$ cherch\'ee, soit il existe un polyn\^ome-cl\'e $\phi$ pour la valuation $\mu _1$ dans $L[x]$ de degr\'e strictement sup\'erieur \` a $1$.  

Si l'ensemble $\Lambda _{\mu}$ n'a pas de plus grand \'el\'ement nous trouvons un sous-ensemble $\{ \delta _{\alpha} ; \alpha \in A \}$ cofinal dans $\Lambda _{\mu}$, index\'e par un ensemble totalement ordonn\'e $A$, sans plus grand \'el\'ement, avec $\delta _{\alpha} < \delta _{\beta}$ pour $\alpha < \beta$, et pour tout $\alpha \in A$ nous choisissons un polyn\^ome $\phi _{\alpha} = x - a_{\alpha}$ v\'erifiant $\mu (\phi _{\alpha}) = \delta _{\alpha}$. 
Alors la famille $\cC$ de valuation d\'efinie par $\cC = \bigl ( \omega _{(a_{\alpha} , \delta _{\alpha})} \bigr ) _{\alpha \in A}$ est une famille continue de valuations de $L[x]$. 
Pour tout $\alpha < \beta$ dans $A$ nous avons $\nu( a_{\beta} - a_{\alpha})= \gamma _{\alpha}$, en particulier la famille $(a_{\alpha})_{\alpha \in A}$ v\'erifie 
$$\nu (a_{\tau} -a_{\sigma}) \ > \ \nu (a_{\sigma} -a_{\rho}) \ $$
pour tous $\tau > \sigma > \rho$ dans $A$, et nous retrouvons la d\'efinition d'Ostrowski de \emph{famille pseudo-convergente} (\cite{Os}, \cite{Ka}).   
De plus si cette famille admet un polyn\^ome-cl\'e limite celui-ci est de degr\'e strictement sup\'erieur \`a $1$.   

Nous en d\'eduisons le r\'esultat suivant. 

\begin{proposition}\label{prop:corpsalgclos}  
Supposons que le corps $K$ est alg\'ebriquement clos, alors une valuation $\mu$ de $K[x]$ est soit une valuation de la forme $\mu = \omega _{(a,\delta )}$, soit une une valuation associ\'ee \`a une famille pseudo-convergente. 
\end{proposition}  

\begin{preuve}
Comme il ne peut pas exister de polyn\^ome-cl\'e ou polyn\^ome-cl\'e limite de degr\'e strictement plus grand que $1$, soit l'ensemble $\Lambda _{\mu} = \left\{ \mu (x-b) \ | \ b \in K \right\}$ a un plus grand \'el\'ement $\delta$, la valuation $\mu$ est bien sp\'ecifi\'ee et est de la forme $\mu = \omega _{(a,\delta )}$, soit l'ensemble $\Lambda _{\mu}$ n'a pas de plus grand \'el\'ement, la valuation $\mu$ n'est pas bien sp\'ecifi\'ee et elle est associ\'ee \`a la famille pseudo-convergente $(a_{\alpha})_{\alpha \in A}$. 

\hfill$\Box$ 
\end{preuve}

En particulier si la valuation $\mu$ n'est pas bien sp\'ecifi\'ee, le corps r\'esiduel $\kappa _{\mu}$ de la valuation $\mu$ est \'egal au corps r\'esiduel $\kappa _{\nu}$ de la valuation $\nu$ et le groupe des ordres $\Gamma _{\mu}$ est \'egal au groupe des ordres $\Gamma _{\nu}$, et nous en d\'eduisons que l'extension de corps valu\'es $(K(x), \mu ) / (K,\nu )$ est \emph{imm\'ediate} (cf. \cite{Ka}).    
Nous \'etudions dans l'{\scshape {\bf Annexe A}} le lien entre les r\'esultats de Kaplansky sur les extensions imm\'ediates et la pr\'esentation des extensions de valuations \`a partir des familles admissibles. 

\vskip .2cm 

Dans le cas o\`u le corps $K$ est alg\'ebriquement clos, pour d\'ecrire toutes les valuations de $K[x]$ prolongeant la valuation $\nu$ de $K$, nous munissons $K$ de la distance ultram\'etrique associ\'ee \`a $\nu$. 
Pour tout $a\in K$ et tout $\delta \in \Gamma$ nous d\'efinissons les boules ferm\'ee $B$ et ouverte $B ^{\circ}$ de centre $a$ et de rayon $\delta$ respectivement par 
$$\begin{array}{c} 
B = B(a,\delta) = \bigl \{ c\in K \ / \ \nu (a-c) \geq \delta \bigr \} \ , \\
B ^{\circ} = B ^{\circ} (a,\delta) = \bigl \{ c\in K \ / \ \nu (a-c) > \delta \bigr \} \ .
\end{array} $$
Comme la distance d\'efinie par la valuation $\nu$ est ultram\'etrique tout \'el\'ement appartenant \`a une boule ouverte ou ferm\'ee est son centre, plus pr\'ecis\'ement si $b\in B ^{\circ} (a,\delta)$, resp. $b\in B (a,\delta)$, alors $B ^{\circ} (a,\delta) = B ^{\circ} (b,\delta)$, resp. $B (a,\delta) = B (b,\delta)$. 

\begin{proposition}\label{prop:boules} 
Toute boule ferm\'ee $B = B(a,\delta)$ d\'efinit une valuation bien sp\'ecifi\'ee $\mu$ de $K[x]$ par $\mu = \omega _{(a,\delta )}$, qui ne d\'epend pas du centre $a$, et toute valuation bien sp\'ecifi\'ee $\mu$ est de cette forme. 

\`A toute valuation $\mu$ qui n'est pas bien sp\'ecifi\'ee, on peut associer une famille d\'ecroissante $(B_{\alpha} = B(a_{\alpha} , \delta _{\alpha})  )_{\alpha \in A}$ de boules ferm\'ees, o\`u $A$ est un ensemble totalement ordonn\'e sans plus grand \'el\'ement, dont l'intersection $\bigcap _{\alpha} B_{\alpha}$ est vide, telle que la valuation $\mu$ est d\'efinie par 
$$ \mu (x-c) \ = \ Sup \left ( \nu ( c-a_{\alpha}) \ ; \ \alpha \in A \right) \ .$$   
\end{proposition} 

\begin{preuve}
La premi\`ere partie concernant les valuations bien sp\'ecifi\'ees $\mu$ est une cons\'equence directe de ce qui pr\'ec\`ede. 

Pour une valuation qui n'est pas bien sp\'ecifi\'ee $\mu$ il reste \`a v\'erifier que l'intersection $\bigcap _{\alpha} B_{\alpha}$ est vide. En effet la valuation $\mu$ est d\'efinie par une famille continue $\cC = \bigl ( \omega _{(a_{\alpha} , \delta _{\alpha})} \bigr ) _{\alpha \in A}$, et par hypoth\`ese pour tout $c\in K$ il existe $\alpha \in A$ tel que $\mu (x-c) < \mu (x -a_{\alpha})=\delta _{\alpha}$, d'o\`u $\nu (a_{\alpha} -c) < \delta _{\alpha}$. 

\hfill$\Box$ 
\end{preuve}

\begin{remark}
Si $\mu _1$ et $\mu _2$ sont deux valuations bien sp\'ecifi\'ees de $K[x]$ d\'efinies respectivement par $\mu _1 = \omega _{(a_1 , \delta _1 )}$ et $\mu _2 = \omega _{(a_2 , \delta _2 )}$, nous avons $\mu _1 \leq \mu _2$ si et seulement si $\delta _1 \leq \delta _2$ et $\delta _1 \leq \nu (a_1 - a_2)$, c'est-\`a-dire si et seulement si $B(a_2, \delta _2) \subset B(a_1, \delta _1) $. 
\end{remark}

\vskip .2cm  

Soit $\mu$ une valuation bien sp\'ecifi\'ee de $K[x]$ de la forme $\mu = \omega _{(a,\delta )}$, sans supposer que le corps $K$ soit alg\'ebriquement clos. Si $\delta$ n'appartient pas au groupe $\Gamma _{\nu}$, le groupe des ordres $\Gamma _{\mu}$ est \'egal \`a $\Gamma _{\nu} \oplus \delta \Z$ et nous avons
$$ rang. \Gamma _{\mu} / \Gamma _{\nu} = rang. rat. \Gamma _{\mu} / \Gamma _{\nu} =1 \quad\hbox{et}\quad 
 \kappa _{\mu}  = \kappa _{\nu} \ .$$
Si $\delta$ appartient au groupe $\Gamma _{\nu}$, alors le groupe des ordres $\Gamma _{\mu}$ est \'egal \`a $\Gamma _{\nu} $ et le corps r\'esiduel $ \kappa _{\mu} $ est une extension transcendante de $\kappa _{\nu}$ engendr\'e par l'image de $b(x-a)$ o\`u $b\in K$ avec $\nu (b) = - \delta$, d'o\`u: 
$$ \Gamma _{\mu} = \Gamma _{\nu} \quad\hbox{et}\quad dim.alg. _{\kappa _{\nu}} \kappa _{\mu} =1\ .$$
Dans tous les cas l'alg\`ebre gradu\'ee associ\'ee ${\rm gr}_{\mu}K[x]$ associ\'ee \`a la valuation bien sp\'ecifi\'ee $\mu = \omega _{(a,\delta )}$ est isomorphe \`a ${G^{(0)}} [S]$, 
o\`u ${G^{(0)}}$ est une alg\`ebre gradu\'ee simple isomorphe \`a ${\rm gr}_{\nu}K$, et $S$ est l'image $H_{\mu}(x-a)$. 
Pour tout $a'\in K$ nous avons 

$$\begin{array}{ll} 
H_{\mu}(x-a') = S & \hbox{si}\ \mu (x-a') = \mu (x-a) < \nu (a-a'), \\
H_{\mu}(x-a') = S - H_{\nu}(a-a') & \hbox{si}\ \mu (x-a') = \mu (x-a) = \nu (a-a'), \\
H_{\mu}(x-a') = H_{\nu}(a-a') & \hbox{si}\ \mu (x-a') = \nu (a-a') < \mu (x-a), 
\end{array} $$

\noindent d'o\`u la remarque suivante. 

\begin{remark} \label{rmq:inversible}
Soit $\mu$ une valuation bien sp\'ecifi\'ee de $K[x]$ de la forme $\mu = \omega _{(a,\delta )}$, alors le polyn\^ome $(x-b)$ a son image  $H_{\mu}(x-b)$ inversible dans l'alg\`ebre gradu\'ee ${\rm gr}_{\mu}K$ si et seulement si $\nu (a-b) < \delta$.
\end{remark} 

\vskip .2cm

      \section{Passage \` a la cl\^oture alg\'ebrique}    
%

Dans cette partie nous nous donnons un corps valu\'e $(K,\nu )$, une cl\^oture alg\'ebrique $\bar K$ de $K$ et une valuation $\bar\nu$ de $\bar K$ qui prolonge la valuation $\nu$, nous consid\'erons une valuation $\mu$ de $K[x]$ qui prolonge la valuation $\nu$ et nous allons \'etudier les prolongements $\bar\mu$ \`a $\bar K[x]$ de $\mu$ qui sont aussi des prolongements de la valuation $\bar \nu$.

\begin{definition} \label{def:envergure1} 
Pour tout polyn\^ome $f$ de $K[x]$ nous d\'efinissons l'\emph{\'etendue} de $f$, que nous notons $\varepsilon _{\mu}(f)$,  par 
$$\varepsilon _{\mu}(f) \ = \ Sup \left ( \bar \mu ( x- a ) \ ;  \hbox{ $a$ racine de $f$ dans $\bar K$} \right ) \ ,$$
o\`u $\bar\mu$ est un prolongement de la valuation $\mu$. 
\end{definition} 

\begin{proposition}\label{prop:independance}
L'\'etendue du polyn\^ome $f$ est ind\'ependante du prolongement $\bar\mu$ choisi. 
\end{proposition} 

\begin{preuve} 
Soient $\bar\mu _1$ et $\bar\mu _2$ deux prolongements de la valuation $\mu$ de $K[x]$ \`a $\bar K [x]$, alors il existe $\sigma$ dans $Aut(\bar K(x)/K(x)) \simeq Aut(\bar K/K)$ tel que $\mu _1 = \mu _2 \circ \sigma$, en particulier pour tout $a \in \bar K$ nous avons 
$$\bar\mu _1 ( x-a) \ = \ \bar\mu _2 ( x- \sigma (a)) \ .$$
Comme le groupe $Aut(\bar K / K)$ agit transitivement sur les racines du polyn\^ome $f$ nous en d\'eduisons que les ensembles $\bigl \lbrace \bar\mu _1 (x-a) \ ;  \hbox{ $a$ racine de $f$ dans $\bar K$} \bigr \rbrace$ et $\bigl \lbrace \bar\mu _2 (x-a) \ ;  \hbox{ $a$ racine de $f$ dans $\bar K$} \bigr \rbrace$ sont \'egaux.  

\hfill \tf   
\end{preuve}

\begin{proposition}\label{prop:equiv} 
La valuation $\bar\mu$ est bien sp\'ecifi\'ee si et seulement si la valuation $\mu$ l'est. 
\end{proposition} 

\begin{preuve} 
Comme la valuation $\bar\mu$ est une extension de la valuation $\mu$ nous avons un morphisme injectif canonique d'alg\`ebres gradu\'ees int\`egres
$$\xymatrix{ \rho : {\rm gr}_{\mu} (K[x]) \ar@{^{(}->}[r] &  {\rm gr}_{\bar\mu} (\bar K[x]) \ .}$$

Nous supposons d'abord que la valuation $\mu$ n'est pas bien sp\'ecifi\'ee, alors l'alg\` ebre gradu\'ee ${\rm gr}_{\mu} (K[x])$ est simple, et nous devons montrer que pour tout $a \in \bar K$, l'\'el\'ement $H_{\bar \mu}(x-a)$ est inversible dans ${\rm gr}_{\bar\mu} (\bar K[x])$. 
Soit $a \in \bar K$, et soit $\phi$ le polyn\^ome minimal de $a$ sur $K$, nous \'ecrivons $\phi (x) = \prod _{i=1 }^d (x-a_i)$, o\`u les $a_i$ sont les racines de $\phi$ dans $\bar K$. 
L'image de $H_{\mu}(\phi )$ par $\rho$ est \'egale au produit 
$$\prod _{i=1 }^d H_{\bar \mu}(x-a_i) \ ,$$  
et comme $H_{\mu}(\phi )$ est inversible dans ${\rm gr}_{\mu} (K[x])$, chacun des facteurs $H_{\bar \mu}(x-a_i)$ est inversible dans ${\rm gr}_{\bar\mu} (\bar K[x])$. 
\vskip .2cm 

R\'eciproquement nous supposons que la valuation $\mu$ est bien sp\'ecifi\'ee, alors l'alg\`ebre gradu\'ee ${\rm gr}_{\mu} (K[x])$ est isomorphe \`a une alg\`ebre de polyn\^omes $G^{(0)}[T]$, o\`u $G^{(0)}$ est une alg\`ebre simple et $T$ est l'image $H_{\mu}(\phi )$ du polyn\^ome-cl\'e $\phi$ qui d\'efinit la valuation $\mu$.

Si la valuation $\bar\mu$ n'\'etait pas bien sp\'ecifi\'ee, l'image de $T$ par $\rho$ serait inversible dans ${\rm gr}_{\bar\mu} (\bar K[x])$ et il existerait en particulier $b \in \bar K$ tel $\rho (T) = H_{\bar\mu} (b) $. 
Il existe un polyn\^ome \`a coefficients dans $K$, dont $b$ est une racine, $c_nX^n + c_{n-1} X^{n-1} + \ldots + c_1 X + c_0$, avec $c_j \in K$ et $c_0 \not= 0$.  
En particulier, si nous nous restreignons aux termes de valuation minimale, il existe $k$, $k\geq 2$, et $0 \leq i_1 < i_2 < \ldots < i_k \leq n$ tels que 
$$H_{\bar\mu} (c_{i_k} b^ {i_k}) + H_{\bar\mu} (c_{i_{k-1}} b^ {i_{k-1}}) + \ldots + H_{\bar\mu} (c_{i_1} b^ {i_1}) \ = \ 0$$ 
dans ${\rm gr}_{\bar\mu} (\bar K[x])$.  
Nous en d\'eduirions la relation 
$$H_{\mu} (c_{i_k}) T^ {i_k} + H_{\mu} (c_{i_{k-1}}) T^ {i_{k-1}} + \ldots + H_{\mu} (c_{i_1}) T^ {i_1} \ = \ 0$$ 
dans ${\rm gr}_{\mu} (K[x]) = G^{(0)}[T]$, ce qui est impossible car les $H_{\mu} (c_{i})$ sont dans $G^{(0)}$. 

\hfill \tf    
\end{preuve}   

\begin{proposition}\label{prop:element-inversible}
Soit $f$ un polyn\^ome dans $K[x]$ alors son image $H_{\mu}(f)$ est inversible dans ${\rm gr}_{\mu} (K[x])$ si et seulement si son image $H_{\bar\mu}(f) = \rho (H_{\mu}(f))$ est inversible dans ${\rm gr}_{\mu} (\bar K[x])$. 
\end{proposition}

\begin{preuve}
L'implication $H_{\mu}(f)$ inversible $\Rightarrow \ H_{\bar\mu}(f)$  inversible est \'evidente. 

Pour montrer l'implication r\'eciproque, nous pouvons supposer que les valuations $\mu$ et $\bar\mu$ sont bien sp\'ecifi\'ees, et soit $f\in K[x]$ tel que $H_{\bar\mu}(f)$  est inversible. 
Alors il existe $b\in \bar K$ tel que $H_{\bar\mu}(f)=H_{\bar\mu}(b)$, et comme dans la d\'emonstration pr\'ec\'edente nous pouvons trouver des entiers $k$, $k\geq 2$, et $0 \leq i_1 < i_2 < \ldots < i_k \leq n$, et des \'el\'ements  $c_{i_1}, \ldots , c_{i_k}$ dans $K$ non nuls tels que 
$$H_{\bar\mu} (c_{i_k} b^ {i_k}) + H_{\bar\mu} (c_{i_{k-1}} b^ {i_{k-1}}) + \ldots + H_{\bar\mu} (c_{i_1} b^ {i_1}) \ = \ 0$$ 
dans ${\rm gr}_{\bar\mu} (\bar K[x])$.  
Nous en d\'eduisons que nous avons 
$$\bigl ( H_{\mu} (c_{i_k} f^ {i_k - i_1-1}) + \ldots + H_{\mu} (c_{i_2} f^ {i_2 - i_1-1}) \bigr ) H_{\mu}(f) + H_{\mu} (c_{i_1})  \ = \ 0 \ ,$$ 
et par cons\'equent $H_{\mu}(f)$ est inversible. 

\hfill \tf    
\end{preuve}

\begin{theorem} 
La valuation $\mu$ est bien sp\'ecifi\'ee  si et seulement si l'ensemble 
$$E_{\mu} = \{ \varepsilon _{\mu} (f) \ | \ f \in K[x] \}$$  
admet un plus grand \'el\'ement $\overline{\varepsilon _{\mu}}$. 

Dans ce cas un polyn\^ome $\phi$ d\'efinit la valuation $\mu$ si et seulement si $\phi$ est un polyn\^ome unitaire de $K[x]$ de degr\'e minimal v\'erifiant $\varepsilon _{\mu} (\phi ) = \overline{\varepsilon _{\mu}}$. 
\end{theorem} 

\begin{preuve} 
D'apr\`es la proposition \ref{prop:equiv} la valuation $\mu$ est bien sp\'ecifi\'ee si et seulement si la valuation $\bar\mu$ l'est, c'est-\`a-dire si et seulement si l'ensemble $\Lambda _{\bar\mu} = \{ \bar\mu (x-b) \ | \ b \in \bar K \}$ admet un plus grand \'el\'ement $\lambda _{\bar\mu}$ (cf. proposition \ref{prop:corpsalgclos}).  

Supposons que la valuation $\mu$ est bien sp\'ecifi\'ee et soit $a \in \bar K$ tel que $\lambda _{\bar\mu} = \bar\mu (x-a)$, alors par d\'efinition pour tout polyn\^ome $f$ dans $K[x]$ nous avons l'in\'egalit\'e $\varepsilon _{\mu} (f) \leq \lambda _{\bar\mu}$.   
Si $\phi$ est le polyn\^ome minimal de $a$ dans $K[x]$, plus g\'en\'eralement si $\phi$ est un polyn\^ome dans $K[x]$ qui admet $a$ comme racine, nous avons par d\'efinition $\bar\mu (x-a) \leq  \varepsilon _{\mu} (\phi )$. 
Nous en d\'eduisons que l'ensemble $E_{\mu}$ admet un plus grand \'el\'ement $\overline{\varepsilon _{\mu}}$ \'egal \`a $\lambda _{\bar\mu}$, et que la valeur $\overline{\varepsilon _{\mu}}$ est atteinte pour tout polyn\^ome $\phi$ de $K[x]$ admettant $a$ comme racine. 

\vskip .2cm 

R\'eciproquement supposons que la valuation $\mu$ n'est pas bien sp\'ecifi\'ee, l'ensemble $\Lambda _{\bar\mu}$ n'admet pas de plus grand \'el\'ement, en particulier pour tout polyn\^ome $f$ de $K[x]$ nous pouvons trouver $a \in \bar K$ tel que $\bar\mu (x-a) > \bar\mu (x-b)$ pour toute racine $b$ de $f$. Un polyn\^ome $g$ dans $K[x]$ admettant $a$ comme racine v\'erifie alors $\varepsilon _{\mu}(g) > \varepsilon _{\mu}(f)$, par cons\'equent l'ensemble $E_{\mu}$ n'admet pas de plus grand \'el\'ement.

\vskip .2cm 

Soit $\phi$ un polyn\^ome de $K[x]$ qui d\'efinit la valuation $\mu$ alors nous avons un isomorphisme  
$${\rm gr}_{\mu}K[x] \ = \ {G^{(0)}} [T]\ ,$$
o\`u ${G^{(0)}}$ est une alg\`ebre gradu\'ee simple et $T$ est l'image $H_{\mu}(\phi )$ du polyn\^ome $\phi$ et nous \'ecrivons comme pr\'ec\'edemment 
$\phi (x) \ = \ \prod _{i=1 }^d (x-a_i)$, 
o\`u les racines $a_i$ de $\phi$ dans $\bar K$ v\'erifient 
$$\varepsilon _{\mu} (\phi) = \bar\mu (x-a_1) = \ldots = \bar\mu (x-a_c) > \bar\mu (x-a_{c+1}) \geq \ldots \geq \bar\mu (x-a_d) \ . $$ 
Nous ne supposons pas que le polyn\^ome $\phi$ est s\'eparable, par cons\'equent les racines $a_1, \ldots , a_d$ de $\phi$ ne sont pas suppos\'ees distinctes. 
L'image $\rho (T)$ de $T$ dans ${\rm gr}_{\bar\mu} (\bar K[x])$ est alors \'egale au produit 
$$\prod _{i=1 }^d H_{\bar\mu}(x-a_i) \ ,$$ 
et nous d\'eduisons de la d\'emonstration de la proposition \ref{prop:equiv} que $\rho (T)$ n'est pas inversible. En particulier un des facteurs $H_{\bar\mu}(x-a_i)$ est non inversible, par cons\'equent nous avons l'\'egalit\'e $\bar\mu ( x-a_i) = \lambda _{\bar\mu} = Sup \{ \bar\mu (x-b) \ | \ b \in \bar K \}$, d'o\`u   $\varepsilon _{\mu} (\phi ) = \lambda _{\bar\mu}$.
\vskip .2cm 

R\'eciproquement supposons que l'ensemble $E_{\mu} = \{ \varepsilon _{\mu} (f) \ | \ f \in K[x] \}$ admet un plus grand \'el\'ement $\overline{\varepsilon _{\mu}}$ et soit $\phi$ un polyn\^ome unitaire de $K[x]$ v\'erifiant $\varepsilon _{\mu} (\phi ) = \overline{\varepsilon _{\mu}}$.  
Nous \'ecrivons encore $\phi = \prod _{i=1} ^d (x-a_i)$ et si $H_{\mu}( \phi )$ \'etait inversible dans ${\rm gr}_{\mu} (K[x])$ nous d\'eduirions comme pr\'ec\'edemment que chacun des facteurs $H_{\bar\mu}(x-a_i)$ serait inversible dans ${\rm gr}_{\bar\mu} (\bar K[x])$, en particulier il existerait $a \in \bar K$ tel que $\bar\mu (x-a_i) < \bar\mu (x-a)$, ce qui contredit l'hypoth\`ese $\varepsilon _{\mu} (\phi ) = \overline{\varepsilon _{\mu}}$.  
  
\hfill \tf    
\end{preuve} 

\begin{corollary} 
Le polyn\^ome irr\'eductible $\phi$ de $K[x]$ d\'efinit la valuation $\mu$ si et seulement si une racine $a$ de $\phi$ dans $\bar K$ d\'efinit une valuation $\bar\mu$ de $\bar K[x]$ qui prolonge $\mu$. 

Plus pr\'ecis\'ement si nous \'ecrivons 
$\phi (x) = \prod _{i=1 }^d (x-a_i)$, 
avec $a = a_1$ et o\`u les racines $a_i$ de $\phi$ dans $\bar K$ v\'erifient 
$$\varepsilon _{\mu} (\phi) = \bar\mu (x-a_1) = \ldots = \bar\mu (x-a_c) > \bar\mu (x-a_{c+1}) \geq \ldots \geq \bar\mu (x-a_d) \ . $$ 
la valuation $\bar\mu$ est la valuation associ\'ee \`a la paire $(a_i, \delta)$ pour $1\leq i \leq c$ et $\delta = \varepsilon _{\mu} (\phi)$. 
\end{corollary}

\hfill \tf    

\vskip .5cm  

Dans la suite nous appelons $p$ l'exposant caract\'eristique du corps $K$, i.e. $p=1$ si $K$ est de caract\'eristique nulle et $p$ est \'egal \`a la caract\'eristique de $K$ sinon. 

Soit $a$ un \'el\'ement de $\bar K$ et soit $\phi$ son polyn\^ome irr\'eductible sur $K$ alors il existe un polyn\^ome irr\'eductible s\'eparable sur $K$, $\phi _{\rm sep}$, et un entier $n \geq 0$ tel que nous ayons l'\'egalit\'e $\phi (x) = \phi _{\rm sep}( x ^{p^n})$. 
Le degr\'e $d_{\rm s}$ du polyn\^ome $\phi _{\rm sep}$ est par d\'efinition le degr\'e de s\'eparabilit\'e de l'extension $(L/K)$, o\`u $L$ est la sous-extension de $\bar K$ engendr\'ee par $a$, $L=K(a) \simeq K[x] /(\phi )$, et nous avons $d=[L:K] = p^n d_{\rm s}$ et $d_{\rm s} = [L:K]_{\rm sep} = [G:H]$, o\`u nous appelons respectivement $G$ et $H$ les groupes de Galois $Gal(\bar K/K)$ et $Gal(\bar K/L)$. En particulier nous pouvons identifier $H$ au sous-groupe $\{\sigma \in G \ / \ \sigma (a)=a \}$.  
Si nous posons $Rac(\phi )= \{a_1, \ldots , a_{d_{\rm s}} \}$ alors nous avons: 
  
$$\phi (x) \ = \ \phi _{\rm sep}(x ^{p^n}) \ = \ \prod _{i=1} ^{d _{\rm s}} (x-a_i) ^{p^n} \ , $$

\vskip .2cm

Comme pr\'ec\'edemment nous consid\'erons une valuation $\mu$ de l'anneau des polyn\^omes $K[x]$ bien sp\'ecifi\'ee, c'est \`a dire de la forme 
$$\mu = [ \mu _{\sharp} \ ; \ \mu (\phi ) = \gamma ] \ ,$$ 
telle que le polyn\^ome-cl\'e ou le polyn\^ome-cl\'e limite $\phi$ est de degr\'e $d$, $d \geq 2$, et nous posons $\phi (x)= \phi _{\rm sep}(x ^{p^n})$ avec $\phi _{\rm sep}$ s\'eparable de degr\'e $d_{\rm s}$.  

Il existe une racine $a_1$ du polyn\^ome $\phi$ dans $\bar K$ et une valuation $\bar\mu$ de l'anneau $\bar K[x]$ qui prolonge la valuation $\mu$ qui est de la forme $\bar\mu = \omega _{(a_1,\delta )}$ pour une certaine valeur $\delta$, et nous notons $a_j$, $1\leq j \leq d_{\rm s}$, les racines de $\phi$ de telle fa\c con que nous ayons 
$$\bar\nu (a_1-a_j) \geq \bar\nu (a_1-a_{j+1})\quad\hbox{pour}\ 2\leq j \leq d_{\rm s}-1 \ .$$
 
Pour toute valeur $\delta$ dans $\tilde\Gamma$ nous notons $c (\delta )$ le plus petit entier $j$, $1\leq j\leq d_{\rm s}$, tel que nous ayons $\bar\nu (a_1-a_j) \geq \delta$, en particulier si $c(\delta )<n$ nous avons 
$$\bar\nu (a_1-a_{c(\delta )}) \geq \delta > \bar\nu (a_1-a_{c(\delta )+1})\ .$$ 
La valuation $\omega _{(a_1,\delta )}$ v\'erifie alors $\omega _{(a_1,\delta )}= \omega _{(a_j,\delta )}$ pour $1 \leq j \leq c(\delta )$ et nous avons les \'egalit\'es suivantes: 
 $$\begin{array}{l}
\omega _{(a_1,\delta )} (x - a_j) = \delta \ \hbox{pour}\ 1 \leq j \leq c(\delta ) \\  
\omega _{(a_1,\delta )} (x - a_j) = \bar\nu (a_1 - a_j) < \delta\ \hbox{pour} \ c(\delta )+1 \leq j \leq d_{\rm s} \ .
 \end{array}$$
 
 \begin{lemma}\label{lem:egalite}
Avec les notations pr\'ec\'edentes supposons  que la valuation $\omega _{(a_1,\delta )}$ soit un prolongement $\bar\mu$ de la valuation $\mu$ \`a $\bar K[x]$, alors nous avons: 
$$\mu (\phi) \ = \ \gamma  \ = \ p^n \biggl ( c(\delta ) \delta + \sum _{j=c(\delta )+1} ^{d_{\rm s}} \bar\nu (a_1 - a_j) \biggr ) \ \leq \ p^n d_{\rm s} \delta \ = \ {\rm deg}\phi \ \delta \ ,$$
avec l'in\'egalit\'e stricte $\mu (\phi ) < {\rm deg}\phi \ \delta$ seulement dans le cas o\`u $c(\delta ) < d_{\rm s}$, c'est-\`a-dire dans le cas o\`u la valuation $\mu$ admet plusieurs prolongements distincts \`a $\bar K[x]$.  
\end{lemma}

\begin{proposition}\label{prop:unicite-de-delta} 
Pour toute racine $a$ du polyn\^ome-cl\'e $\phi$ il existe au plus une valuation $\bar\mu$ de $\bar K[x]$ qui prolonge les valuations $\bar\nu$ et $\mu$ qui soit de la forme $\bar\mu = \omega _{(a,\delta )}$.  
\end{proposition} 

\begin{preuve} 
Soit $a=a_1$ une racine de $\phi$, nous choisissons les autres racines $a_i$ de fa\c con \`a avoir les in\'egalit\'es $\bar\nu (a_1-a_j) \geq \bar\nu (a_1-a_{j+1})$. 
Supposons que nous ayons deux valeurs $\delta$ et $\delta '$ dans $\tilde\Gamma$ telles que les valuations $\omega _{(a_1,\delta )}$ et $\omega _{(a_1,\delta ')}$ soient des prolongements de la valuation $\mu$, et notons respectivement $c$ et $c'$ les plus petits entiers tels que nous ayons $\bar\nu (a_1-a_c) \geq \delta$ et  $\bar\nu (a_1-a_{c'}) \geq \delta '$. 

D'apr\`es le lemme \ref{lem:egalite} nous avons l'\'egalit\'e: 
$$p^n \biggl ( c \delta + \sum _{j=c+1} ^{d_{\rm s}} \bar\nu (a_1 - a_j) \biggr ) \
= \  p^n \biggl ( c' \delta ' + \sum _{j=c'+1} ^{d_{\rm s}} \bar\nu (a_1 - a_j) \biggr ) \ .$$ 

Si $c=c'$ nous en d\'eduisons l'\'egalit\'e $\delta = \delta '$. 

Si $c '>c$ alors nous avons $\delta > \bar\nu (a_1-a_{c'}) \geq \delta '$, et nous d\'eduisons de l'\'egalit\'e pr\'ec\'edente que nous avons  
$$  c \delta + \sum _{j=c+1} ^{c'} \bar\nu (a_1 - a_j) = c' \delta '  \ ,$$ 
avec $\delta > \delta '$ et $\bar\nu (a_1 - a_j) \geq \delta ' $ pour tout $j \leq c'$, ce qui est impossible. 

\hfill \tf   
\end{preuve}

Nous d\'eduisons de ce qui pr\'ec\`ede que si les racines du polyn\^ome $\phi$ sont tr\`es proches pour la distance d\'efinies par $\bar\nu$ sur $\bar K$ la valuation $\mu$, d\'efinie par le polyn\^ome $\phi$ admet un seul prolongement. 

\begin{proposition}\label{prop:unicite-de-mubar} 
La valuation $\mu$ admet un unique prolongement $\bar\mu$ \`a $\bar K[x]$ si et seulement si nous avons pour tout couple de racines $(a,b)$ de $\phi$: 
$$\bar \nu (a-b) \ \geq \ \mu (\phi ) / {\rm deg}\phi \ . $$
\end{proposition} 

\begin{preuve} 
Supposons que la valuation $\mu$ admette un seul prolongement $\bar \mu = \omega _{(a,\delta )}$, alors nous avons pour toute racine $b$ l'in\'egalit\'e $\bar\nu (a-b) \geq \delta$ avec $\delta =  \mu (\phi ) / {\rm deg}\phi $.  

Supposons maintenant que nous avons l'in\'egalit\'e $\bar \nu (a-b) \geq \mu (\phi ) / {\rm deg}\phi$ pour tout couple $(a,b)$ de racines, et soit $\bar \mu$ un prolongement de $\mu$ de la forme $\bar \mu = \omega _{(a,\delta )}$. 
Nous notons $c$ le nombre de racines $b$ v\'erifiant  $\bar\nu (a-b) \geq \delta $ et $d=p^n d_{\rm s}={\rm deg}\phi$, et nous d\'eduisons alors du lemme \ref{lem:egalite} la relation 
$$ \mu (\phi)  \ = \ p^n \biggl ( c \delta + \sum _{ \{b/ \bar\nu (a-b) < \delta \} }  \bar\nu (a - b) \biggr ) \ \geq \ p^n c \delta + p^n (d_{\rm s}-c)  \mu (\phi ) / d \ ,$$ 
d'o\`u  le r\'esultat. 

\hfill \tf   
\end{preuve}

\vskip .2cm 

Dans la suite de ce paragraphe nous supposerons que le corps valu\'e $(K,\nu )$ est hens\'elien, la valuation $\nu$ admet donc un seul prolongement not\'e $\bar\nu$ au corps $\bar K$, en particulier pour tout automorphisme $\sigma$ dans le groupe de Galois $G=Gal(\bar K / K)$ la valuation $\bar\nu \circ \sigma$ est \'egale \`a la valuation $\bar\nu$. 

\begin{theorem}\label{th:prolongements}
Soit $\mu$ une valuation bien sp\'ecifi\'ee de $K[x]$ d\'efinie par un polyn\^ome $\phi$ de degr\'e $d$, et soit $a$ une racine de $\phi$ dans $\bar K$. 
Alors il existe une valeur $\delta$ dans $\tilde\Gamma$ et des automorphismes $\sigma _1=id , \sigma _2 , \ldots , \sigma _k$, avec $1\leq k \leq d$, dans le groupe de Galois $G=Gal(\bar K /K)$ tels que les prolongements $\bar\mu ^{(l)}$ de la valuation $\mu$ \`a $\bar K [x]$ sont les valuations $\omega _{( \sigma _l(a),\delta )}$, pour $1\leq l \leq k$. 
\end{theorem}

\begin{preuve} 
Si $\bar\mu$ est un prolongement de la valuation $\mu$, alors tous les prolongements de $\mu$ sont de la forme $\bar\mu \circ \sigma$ pour $\sigma$ appartenant au groupe de Galois $Gal(\bar K/K)$. 
Supposons que $\bar\mu$ soit la valuation $\omega _{(a,\delta )}$ pour une racine $a$ de $\phi$, alors la valuation $\bar\mu \circ \sigma$ est la valuation $\omega _{(\sigma ^{-1}(a),\delta )}$. 
En effet la valuation $\bar\mu =\omega _{(a,\delta )}$ est d\'etermin\'ee par l'\'egalit\'e 
$$\bar\mu (x-b) \ =\ Inf \bigl ( \delta , \bar\nu (a-b) \bigr ) \ ,$$ 
o\`u $b$ parcourt $\bar K$. 

Comme $(K,\nu )$ est hens\'elien nous en d\'eduisons que pour tout $b$ nous avons 
$$\bar\mu \circ  \sigma (x-b) \ =\ Inf \bigl ( \delta , \bar\nu (a-\sigma (b)) \bigr ) \ =\ Inf \bigl ( \delta , \bar\nu (\sigma ^{-1}(a)-b) \bigr ) \ .$$ 

Comme le groupe de Galois agit transitivement sur les racines du polyn\^ome $\phi$ nous en d\'eduisons que pour toute racine $a$ de $\phi$ il existe un prolongement $\bar\mu$ de la valuation $\mu$ de la forme $\omega _{(a,\delta )}$, et nous d\'eduisons de la proposition \ref{prop:unicite-de-delta} que la valeur $\delta$ est d\'etermin\'ee uniquement par la valuation $\mu$.  

\hfill \tf   
\end{preuve} 

\vskip .2cm

Les diff\'erents prolongements de la valuation $\mu$ \`a $\bar K[x]$ correspondent aux diff\'erentes boules ferm\'ees $B(a, \delta )$ de $\bar K$ o\`u $a$ parcourt l'ensemble $Rac(\phi )$ des racines distinctes de $\phi$, et deux boules $B(a, \delta )$ et $B(a', \delta )$ sont disjointes si $\bar\nu (a-a') < \delta$ ou \'egales sinon. 

Soit $B$ une boule ferm\'ee de $\bar K$, alors pour tout $\sigma$ dans le groupe de Galois $G$ la boule $\sigma .B$ est \'egale \`a $B$ si pour tout $a$ dans $B$ l'image $\sigma (a)$ appartient \`a $B$ et est disjointe de $B$ sinon, en fait il suffit qu'il existe un \'el\'ement $a$ de $B$ tel que son image appartienne \`a $B$ pour que $\sigma .B$ soit \'egale \`a $B$. 
Nous pouvons d\'efinir $H_B$ par 
$$ H_B \ := \ \bigl \{ \sigma \in G \ / \sigma .B = B \bigr \} \ .$$
Il est facile de v\'erifier que les $H_B$ sont des sous-groupes du groupe de Galois $G$ qui v\'erifient
\begin{enumerate} 
\item $H_B \subset H_{B'}$ pour $B \subset B'$; 
\item $H_{\tau .B} = \tau . H_B . \tau ^{-1}$.
\end{enumerate} 
Il y a une bijection naturelle entre l'ensemble quotient $G/H_B$ et l'ensemble $\bigl \{ B^{(l)} \bigr \}$ des boules ferm\'ees disjointes conjugu\'ees \`a $B$ par l'action du groupe de Galois $G$. 

\vskip .2cm

Si un \'el\'ement $b$ de $\bar K$ appartient \`a $B$, la boule $B$ est la boule ferm\'ee $B(b,\delta )$ et le sous-groupe $H_B$ s'identifie au sous-groupe $H_{(b,\delta )}$ d\'efini par 
$$ H_{(b,\delta )} \ := \ \bigl \{ \sigma \in G \ / \bar\nu (b-\sigma (b) ) \geq \delta \bigr \} \ ,$$
et comme pr\'ec\'edemment nous avons $H_{(b,\delta )} \subset H_{(b,\delta ')}$ pour $\delta \geq \delta '$, et $H_{(\tau (b),\delta )} = \tau . H_{(b,\delta )} . \tau ^{-1}$.
Rappelons que pour tout \'el\'ement $b$ de $\bar K$ nous d\'efinissons la \emph{constante de Krasner} $\Delta _K (b)$ par 
$$\Delta _K (b) \ = \ Sup \bigl ( \bar\nu (b -b') \ ; \ b'\ \hbox{conjugu\'e de}\ b \ , \ b'\not=b \bigr ) \ . $$ 
Alors, pour $\delta > \Delta _K (b)$ le sous-groupe $H_{(b,\delta )}$ est \'egal au groupe de Galois $H=Gal(\bar K/K(b))$. 

\vskip .2cm 

La bijection naturelle entre l'ensemble quotient $G/H$ et l'ensemble $Rac(\psi )$ des racines du polyn\^ome irr\'eductible $\psi =Irr _K(b)$ de $b$ sur $K$ induit une bijection entre $H_B/H$ et l'ensemble $Rac_B(\psi )$ des racines de $\psi$ appartenant \`a la boule ferm\'ee $B=B(b,\delta )$. 
En particulier l'ensemble $G/H_B$ est fini et il existe $\sigma _1, \ldots, \sigma _k$ dans $G$ tels que l'ensemble $\bigl \{ B^{(l)} \bigr \}$ des boules ferm\'ees conjugu\'ees \`a $B$ soit \'egal \`a $\bigl \{ \sigma _l .B \bigr \}$, 
avec $\sigma _1=id$ et $B_{(1)}=B$. 
En particulier les prolongements $\bar\mu ^{(l)}$ de la valuation $\mu$ \`a $\bar K[x ]$ d\'efinis au th\'eor\`eme \ref{th:prolongements} sont les valuations associ\'ees aux boules $B^{(l)}$. 

Nous d\'efinissons la \emph{distance} $\eta_{l,l'}$ entre deux boules distinctes $B_{(l)}$ et $B_{(l')}$ par 
$$\eta_{l,l'} = \bar\nu ( b-b') \ \hbox{pour $b\in B_{(l)}$ et $b'\in B_{(l')}$} \ , $$ 
et ceci est ind\'ependant des \'el\'ements $b$ et $b'$ choisis car la distance associ\'ee \`a $\bar\nu$ est ultra-m\'etrique, et v\'erifie $\eta_{l,l'} < \delta$. 

\vskip .2cm 

Soit $b$ appartenant \`a $\bar K$, nous notons $L$ l'extension $K(b)$, $\psi$ son polyn\^ome irr\'eductible sur $K$, $\psi _{\rm sep}$ le polyn\^ome irr\'eductible s\'eparable associ\'e, et nous posons $d=p^n d_{\rm s}$ o\`u $d={\rm deg}\ \psi =[L:K]$ et $d_{\rm s}={\rm deg}\ \psi_{\rm sep} = [L:K]_{\rm sep} = [G:H]$, avec  $H=Gal(\bar K/L)$.   
 
\begin{proposition}\label{prop:valeur-de-mu}
Si $b$ appartient \`a la boule $B$ nous avons l'inclusion $H\subset H_B$ et si nous posons $c=[H_B :H]$ et $k=[G: H_B]$, nous avons l'\'egalit\'e 
$$d_{\rm s} = kc \ .$$ 
De plus nous pouvons trouver $\delta _1, \delta _2, \ldots , \delta _c$ dans $G$ avec $\delta _1=id$, tels que pour tout $l$, $1\leq l\leq k$, l'ensemble $Rac_{B^{(l)}}(\phi )$ des racines de $\phi$ appartenant \`a la boule $B^{(l)}$ est \'egal \`a l'ensemble $\bigl\{ \sigma _l\delta _i(b) , 1 \leq i \leq c \bigr \}$.   
\end{proposition} 

\begin{preuve} 
C'est une cons\'equence de l'\'egalit\'e $H_B= H_{(b,\delta )}$ et du fait que la valuation $\bar\nu\circ\sigma$ est \'egale \`a la valuation $\bar\nu$ pour tout $\sigma$ dans $G$. 

\hfill \tf   
\end{preuve} 
\vskip .2cm 

\begin{corollary} \label{cor:valeur-de-mu}
Soient $\mu$ une valuation bien sp\'ecifi\'ee de $K[x]$, $\bar\mu$ un prolongement de $\mu$ \`a $\bar K[x]$ associ\'e \`a la boule ferm\'ee $B$ de diam\`etre $\delta$ et $b$ un \'el\'ement de $\bar K$ appartenant \`a $B$ de polyn\^ome irr\'eductible sur $K$ $Irr_K(b) = \psi$. Alors avec les notations pr\'ec\'edentes nous avons 
$$\mu(\psi ) \ = \ p^n c \biggl ( \delta + \sum _{l=2} ^k \eta_{1,l} \biggr ) \ \leq \ {\rm deg} \psi . \delta \ . $$
\end{corollary} 

\hfill \tf   
\vskip .2cm 

En particulier nous retrouvons que si la valuation $\mu$ est d\'efinie par le polyn\^ome $\phi$ nous avons l'in\'egalit\'e $\gamma = \mu (\phi )\leq {\rm deg} \phi . \delta$, avec \'egalit\'e si et seulement si la valuation $\mu$ admet un seul prolongement \`a $\bar K[x]$. 

\vskip .2cm

D'apr\`es la proposition \ref{prop:boules} nous pouvons associer \`a une valuation de la forme $\bar\mu = \omega _{(a,\delta )}$ la boule ferm\'ee $B= B(a,\delta )$ de $\bar K$, et la valuation est enti\`erement d\'etermin\'ee par $B$, 
et nous d\'eduisons du th\'eor\`eme \ref{th:prolongements} nous pouvons associer \`a la valuation $\mu$ de $K[x]$ une famille ${\mathcal B}(\mu ) = \bigl ( B^{(l)} \bigr ) _{1\leq l\leq k} $ de boules ferm\'ees disjointes, de m\^eme diam\`etre $\delta$, chacune des boules $B^{(l)}$ correspondant \`a la valuation $\bar\mu ^{(l)}$. 

\begin{proposition} \label{pr:noinversible-boule} 
Soit $\mu$ une valuation bien sp\'ecifi\'ee de $K[x]$ et soit ${\mathcal B}(\mu ) = \bigl ( B^{(l)} \bigr ) _{1\leq l\leq k} $ la famille de boules ferm\'ees de $\bar K$ associ\'ee. 
Alors un polyn\^ome $f$ de $K[x]$ a son image $H_{\mu}(f)$ non inversible dans ${\rm gr}_{\mu} (K[x])$ si et seulement si l'ensemble $Rac(f)$ des racines de $f$ est inclus dans la r\'eunion des boules ferm\'ees $\bigcup _{1 \leq l \leq k} B^{(l)}$. 
\end{proposition}
 
\begin{preuve} 
Soit $\bar\mu$ un prolongement de $\mu$ \`a $\bar K [x]$ de la forme $\bar\mu = \omega _{(a,\delta )}$, d'apr\`es la proposition \ref{prop:element-inversible} le polyn\^ome $f$ a son image $H_{\mu}(f)$ non inversible dans ${\rm gr}_{\mu} (K[x])$ si et seulement si $f$ a une racine $b$ tel que $\bar\nu (a-b) \geq \delta$. 
Comme le groupe de Galois agit transitivement sur l'ensemble des racines $Rac(f)$ et sur l'ensemble des boules ${\mathcal B}(\mu )$, nous en d\'eduisons le r\'esultat. 

\hfill \tf   
\end{preuve} 

\vskip .2cm 

\begin{remark}\label{rmq:divisibilite}  
A toute valuation bien sp\'ecifi\'ee $\mu$ de $K[x]$, nous pouvons associer l'entier $k$ d\'efini comme le nombre de boules ferm\'ees distinctes $B^{(l)}$ conjugu\'ees \`a la boule ferm\'ee $B$ associ\'ee \`a un prolongement $\bar\mu$ de $\mu$ \`a $\bar K[x]$. 

Nous d\'eduisons de ce qui pr\'ec\`ede que tout polyn\^ome irr\'eductible $\psi$  de $K[x]$ ayant une racine dans $B$ est de degr\'e divisible par $k$. 
En particulier si $k$ ne divise pas le degr\'e d'un polyn\^ome $f$ son image $H_{\mu}(f)$ est inversible dans ${\rm gr}_{\mu} (K[x])$.  
\end{remark} 

\vskip .2cm 

Nous voulons \'etudier les prolongements $\bar\mu _i$ des valuations $\mu _i$ appartenant \`a une famille admise ${\mathcal A} = \bigl ( \mu _i \bigr ) _{i \in I}$ de valuations de $K[x]$. 
Pour cela il nous faut d'abord \'etudier les prolongements \`a $\bar K[x]$ de deux valuations dont l'une est valuation augment\'ee de l'autre. 

\begin{theorem} \label{th:prolongement-valuation-augmentee} 
Soit $\mu$ une valuation bien sp\'ecifi\'ee de $K[x]$ obtenue comme valuation augment\'ee $\mu = [ \mu _{\sharp} \ ; \ \mu (\phi ) = \gamma ]$ et soit $\bar \mu _{\sharp}$ un prolongement de la valuation $\mu _{\sharp}$ \`a $\bar K [x]$. 
Alors il existe un prolongement $\bar\mu$ de la valuation $\mu$ \`a $\bar K [x]$ qui est obtenu comme valuation augment\'ee $\bar\mu = [ \bar\mu _{\sharp} \ ; \ \bar\mu (\bar\phi ) = \bar\gamma ]$. 
\end{theorem}

\begin{preuve} 
Comme la valuation $\mu _{\sharp}$ admet une valuation augment\'ee elle est bien sp\'ecifi\'ee, alors pour tout prolongement $\bar\mu _{\sharp}$ il existe une racine $a _{\sharp}$ du polyn\^ome-cl\'e $\phi _{\sharp}$ d\'efinissant la valuation $\mu _{\sharp}$ et une valeur $\delta _{\sharp}$ dans $\tilde\Gamma$ telles que le prolongement $\bar\mu _{\sharp}$ soit la valuation $\omega _{(a _{\sharp} , \delta _{\sharp})}$. 

Comme le polyn\^ome $\phi$ est un polyn\^ome-cl\'e pour la valuation $\mu _{\sharp}$ son image $H_{\mu _{\sharp}}(\phi )$ dans l'anneau gradu\'e  ${\rm gr} _{\mu _{\sharp}} K[x]$ n'est pas inversible. 
Par cons\'equent il existe au moins une racine $a$ de $\phi$ dans $\bar K$ telle que $H_{\bar\mu _{\sharp}}(x-a )$ ne soit pas inversible dans l'anneau gradu\'e  ${\rm gr} _{\bar\mu _{\sharp}} \bar K[x]$, d'o\`u d'apr\`es la remarque \ref{rmq:inversible} l'in\'egalit\'e $\bar\nu (a - a _{\sharp}) \geq \delta _{\sharp}$. 
Nous en d\'eduisons que la valuation $\bar\mu _{\sharp}$ peut s\'ecrire $\bar\mu _{\sharp} = \omega _{(a, \delta _{\sharp})}$

Soit $\bar\mu$ le prolongement de la valuation $\mu$ associ\'e \`a la racine $a$ du polyn\^ome $\phi$, c'est-\`a-dire qui est de la forme $\bar\mu = \omega _{(a,\delta )}$. 
Deux valuations de la forme $\rho _1 = \omega _{(a,\delta _1)}$ et $\rho _2 = \omega _{(a,\delta _2)}$ sont comparables et v\'erifient $\rho _1 \leq \rho _2$ si et seulement si $\delta _1 \leq \delta _2$. Par cons\'equent comme nous avons l'in\'egalit\'e $\mu _{\sharp} \leq \mu$ nous avons les in\'egalit\'es $\delta _{\sharp} \leq \delta$ et $\bar\mu _{\sharp} \leq \bar\mu$. 
Nous en d\'eduisons le r\'esultat car les seuls polyn\^omes-cl\'es sur $\bar K[x]$ sont des polyn\^omes de degr\'e $1$, et nous pouvons prendre $\bar\phi = x-a$. 

\hfill \tf   
\end{preuve} 

\vskip .2cm  

\begin{proposition}
Soient $\bar\mu _{\sharp}=  \omega _{(a _{\sharp},\delta _{\sharp})}$ et $\bar\mu =  \omega _{(a,\delta )}$ les deux valuations d\'efinies dans le th\'eor\`eme pr\'ec\'edent, et supposons que les polyn\^omes $\phi _{\sharp}$ et $\phi$ ne sont pas $\mu _{\sharp}$-\'equivalents, alors nous avons l'\'egalit\'e: 
$$\bar\nu ( a - a_{\sharp}) = \delta _{\sharp}$$
\end{proposition} 

\begin{preuve} 
Nous d\'eduisons de la proposition \ref{prop:element-inversible} que pour toute valuation $\mu$ de $K[x]$ et tout prolongement $\bar\mu =  \omega _{(a,\delta )}$ de $\mu$ \`a $\bar K[x]$, si un polyn\^ome $f$ de $K[x]$ n'est pas $\mu$-inversible alors il existe une racine $b$ de $f$ telle que $\bar\nu (a-b) \geq \delta$. 

Comme les polyn\^omes $\phi _{\sharp}$ et $\phi$ ne sont pas $\mu _{\sharp}$-\'equivalents, le polyn\^ome $\phi _{\sharp}$ est $\mu$-inversible, nous en d\'eduisons l'in\'egalit\'e $\bar\nu (a-a_{\sharp}) < \delta$, par cons\'equent nous avons $\bar\mu (x-a_{\sharp}) = \bar\nu (a-a_{\sharp}) < \delta = \bar\mu (x-a)$. Nous d\'eduisons aussi du fait que $\phi _{\sharp}$ est $\mu$-inversible, l'\'egalit\'e $\delta _{\sharp} = \bar\mu _{\sharp} (x-a_{\sharp}) = \bar\mu (x-a_{\sharp})$, d'o\`u le r\'esultat. 

\hfill \tf   
\end{preuve}

 \vskip .2cm

Soit $\mu$ une valuation bien sp\'ecifi\'ee obtenue comme valuation augment\'ee de la valuation $\mu _{\sharp}$ telle que les polyn\^omes $\phi _{\sharp}$ et $\phi$ ne sont pas $\mu _{\sharp}$-\'equivalents, nous appelons ${\mathcal B}(\mu _{\sharp}) = \bigl ( B^{(l)} _{\sharp}\bigr ) _{1\leq l\leq k _{\sharp}} $ et ${\mathcal B}(\mu ) = \bigl ( B^{(l)} \bigr ) _{1\leq l\leq k} $ les familles de boules ferm\'ees de $\bar K$ associ\'ees respectivement aux valuations $\mu _{\sharp}$ et $\mu$.  

Soit $\bar\mu _{\sharp}^{(r)}$ le prolongement de la valuation  $\mu _{\sharp}$ associ\'e \`a la boule $B^{(r)} _{\sharp} = B( a_{\sharp}^{(r)} ,\delta _{\sharp})$, alors pour toute racine $a^{(l)}$ du polyn\^ome $\phi$ v\'erifiant $\bar\nu ( a^{(l)} - a_{\sharp}^{(r)} ) \geq \delta _{\sharp}$ la boule ferm\'ee associ\'ee $B^{(l)} = B( a^{(l)} ,\delta )$ est incluse dans $B^{(r)} _{\sharp}$. 
Comme le groupe de Galois $G=Gal(\bar K /K)$ agit transitivement sur les ensembles des racines $\bigl \{ a_{\sharp}^{(r)} \bigr \}$ et $\bigl \{ a^{(l)}\bigr \}$ respectivement des polyn\^omes $\phi _{\sharp}$ et $\phi$ nous en d\'eduisons le r\'esultat suivant. 

\begin{proposition}\label{prop:boules-incluses} 
Chaque boule $B_{\sharp}^{(r)}$ de ${\mathcal B} (\mu _{\sharp})$ contient $s$ boules $B^{(l)}$ appartenant \`a ${\mathcal B} (\mu )$, o\`u $s$ est un entier strictement positif ind\'ependant de la boule $B_{\sharp}^{(r)}$ choisie, et toute boule  $B^{(l)}$ de ${\mathcal B} (\mu )$ est contenue dans une boule $B_{\sharp}^{(r)}$ de ${\mathcal B} (\mu _{\sharp})$. 
\end{proposition} 

\hfill \tf

 \vskip .2cm
 
 Nous pouvons repr\'esenter les boules dans $\bar K$ associ\'ees aux prolongements des valuations $\mu _{\sharp}$ et $\mu$ par le diagramme suivant:
 \vskip .2cm
 
$$  
\begin{tikzpicture}[very thick,scale=0.37, baseline,>=stealth]
 
\draw (0,0) circle (3.5);
\draw (1.5,1.5) circle (1);
\draw (1.5,-1.5) circle (1);
\draw (-1.5,1.5) circle (1);
\draw (-1.5,-1.5) circle (1);
\fill (6.7:3.8) circle (0pt) node[left] {{\small$B^{(1)}$}};
\fill (-44:5) circle (0pt) node[left] {{\small$B_{\sharp}^{(1)}$}};

\fill (45:2) circle (2.5pt) node[below] {{$a$}};
\fill (0:0) circle (2.5pt) node[below] {{$a_{\sharp}$}};

\draw (10,0) circle (3.5);
\draw (11.5,1.5) circle (1);
\draw (11.5,-1.5) circle (1);
\draw (8.5,1.5) circle (1);
\draw (8.5,-1.5) circle (1);
\fill (2.2:13.8) circle (0pt) node[left] {{\small$B^{(l)}$}};
\fill (-15:13.8) circle (0pt) node[left] {{\small$B_{\sharp}^{(r)}$}};

\draw (5,-8) circle (3.5);
\draw (6.5,-6.5) circle (1);
\draw (6.5,-9.5) circle (1);
\draw (3.5,-6.5) circle (1);
\draw (3.5,-9.5) circle (1);
\fill (-40.6:11.6) circle (0pt) node[left] {{\small$B^{(q)}$}};
\fill (-56:14.2) circle (0pt) node[left] {{\small$B_{\sharp}^{(p)}$}};

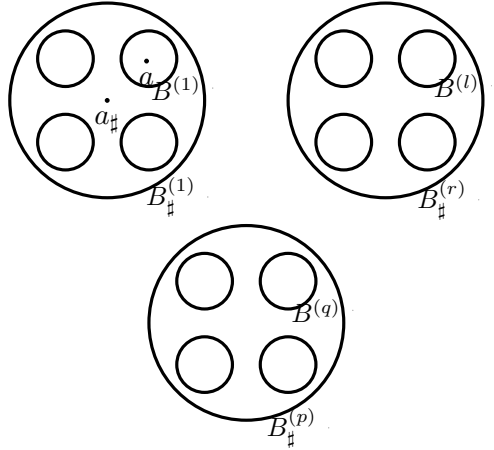
\fill (-70:14.2) circle (0pt) node {{Disposition des boules ferm\'ees des familles ${\mathcal B} (\mu _{\sharp})$ et ${\mathcal B} (\mu )$ dans $\bar K$}};
\end{tikzpicture} 
$$

\vskip .4cm 

Dans le cas o\`u les polyn\^omes $\phi _{\sharp}$ et $\phi$ sont $\mu _{\sharp}$-\'equivalents nous obtenons un r\'esultat identique. 
En effet consid\'erons deux valuations bien sp\'ecifi\'ees $\mu _{\sharp}$ et $\mu$ de $K[x]$ d\'efinies par le m\^eme polyn\^ome $\phi$ et par des valeurs $\gamma _{\sharp}$ et $\gamma$ diff\'erentes avec $\gamma _{\sharp} < \gamma$. 
Soit $a$ une racine de $\phi$ dans $\bar K$ et choisissons deux prolongements $\bar\mu _{\sharp}$ et $\bar\mu$ respectivement de $\mu _{\sharp}$ et $\mu$ \`a $\bar K[x]$, de la forme $\bar\mu _{\sharp}= \omega _{(a, \delta _{\sharp})}$ et  $\bar\mu = \omega _{(a, \delta )}$. 
Comme les valuations $\omega _{(a, \delta _{\sharp})}$ et  $\omega _{(a, \delta )}$ sont comparables et comme nous avons $\mu _{\sharp} \leq \mu$ nous en d\'eduisons l'in\'egalit\'e $\delta _{\sharp} < \delta$. 
Si nous appelons ${\mathcal I} (\mu _{\sharp}) = (I_{\sharp} ^{(l)}\ ; \ 1\leq l \leq k _{\sharp})$ et ${\mathcal I} (\mu )=(I ^{(l)}\ ; \ 1\leq l \leq k)$ les partitions de l'ensemble des racines de $\phi$ associ\'ees respectivement aux valuations $\mu _{\sharp}$ et $\mu$, alors ${\mathcal I} (\mu )$ est plus fine que ${\mathcal I} ( \mu _{\sharp})$, et nous trouvons la m\^eme disposition des boules ferm\'ees associ\'ees aux valuations que pr\'ec\'edemment. 

\vskip .2cm 

Nous voulons \'etudier maintenant le cas d'une famille continue ${\mathcal C} = \bigl ( \mu _{\alpha} \bigr ) _{\alpha \in A}$ de valuations de $K[x]$, 
pour tout $\alpha$ dans $A$ nous notons ${\mathcal B}(\mu _{\alpha}) = \bigl ( B^{(l)} _{\alpha}\bigr ) _{1\leq l\leq k _{\alpha}} $ la famille de boules ferm\'ees de $\bar K$ associ\'ee \`a la valuation $\mu _{\alpha}$ et  $B_{\alpha }$ la r\'eunion  $B_{\alpha} = \bigcup _{1 \leq l \leq k_{\alpha}} B^{(l)} _{\alpha}$. 
Nous d\'eduisons de la proposition \ref{prop:boules-incluses} qu'il existe $\alpha _0$ tel que pour tout $\beta \geq \alpha \geq \alpha _0$ dans $A$ nous avons $k_{\beta} = k _{\alpha} $ et chaque boule $B^{(l)} _{\alpha}$ de ${\mathcal B}(\mu _{\alpha})$ contient une unique boule $B^{(l)} _{\beta}$ de ${\mathcal B}(\mu _{\beta})$, et notons $k _{\mathcal C} = k _{\alpha}$ pour $\alpha \geq \alpha _0$. 
Nous posons alors 
$$B^{(l)}_{\mathcal{C}}  = \bigcap _{\alpha \in A, \alpha \geq \alpha _0} B^{(l)} _{\alpha} \quad \hbox{et} \quad  B_{\mathcal{C}}  = \bigcap _{\alpha \in A} B _{\alpha} = \bigcup _{1 \leq l \leq k _{\mathcal C}} B^{(l)} _{\mathcal{C}} \ .$$   
L'ensemble $B^{(l)}_{\mathcal{C}}$ obtenu comme intersection d'une famille d\'ecroissante de boules ferm\'ees est une boule ferm\'ee, \'eventuellement vide. 
Le groupe de Galois $G = Gal ( \bar K /K )$ agit transitivement sur les racines des polyn\^omes $\phi _{\alpha}$ d\'efinissant les valuations $\mu _{\alpha}$, par cons\'equent toutes les boules $B^{(l)} _{\alpha}$ sont isomorphes.

\begin{proposition}\label{prop:boule-famille-continue} 
Le polyn\^ome $\phi$ appartient \`a $\tilde\Phi \bigl ( \mathcal{C} \bigr ) = \{ f \in K[x] \ | \ \mu _{\alpha} (f) < \mu _{\beta}(f) , \forall \alpha < \beta \in A \}$ si et seulement si l'ensemble des racines $Rac(\phi )$ de $\phi$ est inclus dans $B_{\mathcal C}$.
\end{proposition} 

\begin{preuve} 
Si $\phi$ appartient \`a $\tilde\Phi \bigl ( \mathcal{C} \bigr )$, pour tout $\alpha$ dans $A$ l'image de $\phi$ dans ${\rm gr}_{\mu _{\alpha}} K[x]$ est non inversible et nous d\'eduisons le r\'esultat de la proposition \ref{pr:noinversible-boule}.  

\hfill \tf   
\end{preuve}

 \vskip .2cm

Soient $\mu$ une valuation de $K[x]$ et ${\mathcal A} =  \bigl ( \mu _i \bigr )_{i \in I}$ la famille admise de valuations associ\'ee, chaque valuation $\mu _i$ de la famille $\mathcal A$ est bien sp\'ecifi\'ee, d\'efinie par le polyn\^ome-cl\'e $\phi _i$, et soit  ${\mathcal B} _i = {\mathcal B}(\mu _i)$ la famille de boules ferm\'ees de $\bar K$ asspoci\'ee \`a $\mu _i$ par le th\'eor\`eme \ref{th:prolongements}.

\begin{proposition}\label{prop:familles-incluses} 
La famille $\bigl ( {\mathcal B} _i \bigr )_{i \in I}$ est d\'ecroissante, c'est-\`a-dire pour tout $i<j$ dans $I$ 
chaque boule $B_{i}^{(r)}$ de ${\mathcal B} _i$ contient $s$ boules $B_j^{(l)}$ appartenant \`a ${\mathcal B} _j$, o\`u $s$ est un entier strictement positif ind\'ependant de la boule $B_i^{(r)}$ choisie, et toute boule  $B_j^{(l)}$ de ${\mathcal B} _j$ est contenue dans une boule $B_i^{(r)}$ de ${\mathcal B} _i$. 
\end{proposition} 

\begin{preuve}
C'est une cons\'equence des propositions \ref{prop:boules-incluses} et \ref{prop:boule-famille-continue}.  

\hfill \tf   
\end{preuve}

 \vskip .2cm

Soit $\mu$ une valuation de $K[x]$ et soit ${\mathcal A} =\bigl (\mu _i\bigr )_{i \in I}$ la famille admise associ\'ee. 

\begin{definition} 
Nous appelons \emph{ensemble caract\'eristique} de la valuation $\mu$ le sous-ensemble ${\mathbf B}={\mathbf B}(\mu )$ de $\bar K$ obtenu comme l'intersection des ensembles $B_i$ pour $i$ dans $I$, o\`u chaque $B _i$ est la r\'eunion des boules ferm\'ees $B_{i}^{(r)}$ appartenant \`a ${\mathcal B}(\mu _i)$. 
\end{definition} 

Il existe $i_0 \in I$ tel que les ensembles ${\mathcal B}(\mu _i)$ pour $i \geq i_0$ ont tous le m\^eme nombre $s$ d'\'el\'ements, c'est-\`a-dire que pour tout $j \geq i \geq i_0$ chaque boule ferm\'ee $B_{i}^{(r)}$ de ${\mathcal B}(\mu _i)$ contient une unique boule ferm\'ee $B_{j}^{(r)}$ de ${\mathcal B}(\mu _j)$. 
Pour tout $r$, nous appelons ${\mathcal B}^{(r)}$ la famille d\'ecroissante de boules ferm\'ees $\bigr ( B_{i}^{(r)} \bigl ) _{i \geq i_0}$, et nous notons ${\mathbf B}^{(r)}$ l'intersection 
$${\mathbf B}^{(r)} \ = \ \bigcap _{i \geq i_0} B_{i}^{(r)} \ .$$
Le groupe de Galois $Gal(\bar K/K)$ agit transitivement sur l'ensemble des familles ${\mathcal B}^{(r)}$, donc ces familles sont isomorphes et nous avons l'\'egalit\'e: 
$${\mathbf B} \ = \ \bigcup _{1 \leq r \leq s} {\mathbf B}^{(r)} \ .$$

D'apr\`es la proposition \ref{prop:boules} chaque famille de boules ferm\'ees ${\mathcal B}^{(r)}$ d\'efinit une valuation $\bar \mu ^{(r)}$ de $\bar K[x]$, et nous d\'eduisons de ce qui pr\'ec\`ede le r\'esultat suivant. 

\begin{theorem}\label{th:extensions} 
Les valuations $\bar \mu ^{(r)}$ associ\'ees aux familles ${\mathcal B}^{(r)}$ sont les extensions \`a $\bar K[x]$ de la valuation $\mu$. 

Si l'ensemble caract\'eristique ${\mathbf B}(\mu )$ de $\mu$ est vide la valuation $\mu$ n'est pas bien sp\'ecifi\'ee, et chaque valuation $\bar \mu ^{(r)}$ est d\'efinie par la suite pseudo-convergente associ\'ee \`a la famille ${\mathcal B}^{(r)}$. 

Si l'ensemble caract\'eristique ${\mathbf B}(\mu )$ de $\mu$ est non vide la valuation $\mu$ est bien sp\'ecifi\'ee, et chaque valuation $\bar \mu ^{(r)}$ est d\'efinie par la boule ferm\'ee non vide ${\mathbf B}^{(r)}$. 
\end{theorem}

\begin{corollary}\label{cor:extensions}
Soit $\mu$ une valuation de $K[x]$ et soit ${\mathcal A} = \bigl ( \mu _i \bigr ) _{i \in I}$ une famile admise associ\'ee, alors pour toute extension $\bar\mu$ de $\mu$ \`a $\bar K[x]$, il existe une famille  ${\mathcal A}_{\bar\mu} = \bigl (\bar \mu _i \bigr ) _{i \in I}$ de valuations de $\bar K[x]$ correspondant \`a une suite d\'ecroissante $\bigl ( B _i \bigr ) _{i \in I}$ de boules ferm\'ees de $\bar K$ telle que pour tout $i \in I$ la valuation $\bar\mu _i$ est une extension de la valuation $\mu _i$. 
\end{corollary}

\vskip .2cm

\begin{remark}\label{rmq:carousel}
Soit $(C,0)$ une singularit\'e de courbe plane dans ${\mathbb A}_k^2$ d\'efinie par un polyn\^ome $f \in k[x,y]$ pour un choix judicieux des coordonn\'ees $x$ et $y$, nous consid\'erons alors $f$ comme un \'el\'ement $P$ de $K[x]$ o\`u $K$ est le corps $k(y)$ que nous munissons de la valuation $y$-adique $\nu$. 
Alors l'\'etude de la singularit\'e $(C,0)$ est li\'ee \`a l'\'etude des prolongements $\mu$ de la valuation $\nu$ au corps $L= K[x]/(P)$ et les valuations $\mu _i$ apparaissant dans la famille admise $\cA$ associ\'ee \`a la pseudo-valuation $\bar\mu$ de $K[x]$ de noyau $(P)$ d\'efinie par $\mu$ sont en reli\'ees aux \emph{paires de Puiseux} de la singularit\'e $(C,0)$ (cf. Exemple 3.2 de \cite{Va 4}). 

Il est alors possible de voir une analogie entre la famille des boules $\bigl ( B _i \bigr ) _{i \in I}$ d\'efinie \`a la proposition \ref{prop:familles-incluses} et l'action du groupe de Galois sur celle-ci, et la construction du \emph{carousel} donn\'ee par L\^e  D.T. (voir par exemple \cite{Le 1} ou \cite{Le 2}). 
\end{remark}

\vskip .2cm 

\begin{remark}\label{rmq:independance} 
Soient $\nu$ une valuation de $K$, $\mu$ un prolongement de $\nu$ \`a l'extension $K(x)$ de $K$ et $\cA =\bigl (\mu _i \bigr )_{i \in I}$ la famille de valuations de $K(x)$ associ\'ee \`a $\mu$ consid\'er\'ee comme valuation de l'anneau des polyn\^omes $K[x]$. 
D'apr\`es ce qui pr\'ec\`ede la famille $\cA$ peut \^etre d\'efinie en consid\'erant la famille des boules $\bigl ( B _i \bigr ) _{i \in I}$ de $\bar K$, en particulier on peut en d\'eduire que cette famille ne d\'epend pas du g\'en\'erateur $x$ de l'extension $K(x)$ de $K$. 

Nous retrouvons ainsi le r\'esultat principal de \cite{Va 6}.  
\end{remark}

\section{Restriction d'une valuation d\'efinie sur $\bar K[x]$}

Dans cette partie nous nous donnons une valuation $\bar\mu$ d\'efinie sur $\bar K[x]$ et nous voulons \'etudier sa restriction $\mu$ \`a l'anneau $K[x]$. 
Comme pr\'ec\'edemment nous supposons que $(K,\nu )$ est un corps valu\'e hens\'elien et la valuation $\bar\mu$ de $\bar K[x]$ est un prolongement de l'unique extension $\bar\nu$ de $\nu$ \`a la cloture alg\'ebrique $\bar K$ de $K$. 

Pour tout $b$ appartenant \`a $\bar K$ nous d\'efinissons le degr\'e de $b$ sur $K$ comme le degr\'e de l'extension alg\'ebrique qu'il engendre, c'est aussi la degr\'e du polyn\^ome irr\'eductible ${\rm Irr} _K(b)$ de $b$ sur $K$:  
$${\rm deg}_K (b) = \deg {\rm Irr} _K(b) = [K(b):K] \ .$$

\begin{proposition}\label{prop:degre-minimal} 
Soient $\bar\mu$ une valuation bien sp\'ecifi\'ee de $\bar K[x]$ associ\'ee \`a une boule ferm\'ee $B$, $\mu$ la restriction de $\bar\mu$ \`a $K[x]$ et $\phi$ un polyn\^ome de $K[x]$ d\'efinissant la valuation $\mu$. 
Alors pour tout $b\in B$ nous avons l'in\'egalit\'e:
$${\rm deg}_K (b) \geq \deg \phi \ .$$
\end{proposition}

\begin{preuve}
Nous d\'eduisons de la remarque \ref{rmq:inversible} que si $b$ appartient \`a $B$ l'image de $(x-b)$ est non-inversible dans $ {\rm gr} _{\bar\mu} \bar K[x]$, 
alors l'image du polyn\^ome ${\rm Irr} _K(b)$ est non inversible dans $ {\rm gr} _{\mu} K[x]$, et par cons\'equent il v\'erifie l'in\'egalit\'e $\deg {\rm Irr} _K(b)) \geq \deg \phi$. 

\hfill \tf   
\end{preuve}  
 \vskip .2cm
 
Pour tout sous-ensemble $E$ de $\bar K$ nous d\'efinissons le \emph{degr\'e de $E$ sur $k$} par  
$${\rm deg}_K(E) : = Inf \bigl ({\rm deg}_K (b)\ /\ b \in E \bigr )\ . $$

\begin{theorem}\label{th:polynome-definissant-mu} 
Soient $\bar\mu$ une valuation bien sp\'ecifi\'ee de $\bar K[x]$ et $B$ la boule ferm\'ee de $\bar K$ associ\'ee. 
Alors un polyn\^ome irr\'eductible $\phi$ de $K[x]$ d\'efinit la valuation $\mu$ restriction de $\bar\mu$ \`a $K[x]$ si et seulement si $\phi$ a une racine $a$ appartenant \`a $B$ avec 
$${\rm deg}_K(a) \ = \ {\rm deg}_K(B)\ . $$
\end{theorem} 

\begin{preuve}
Soit $\phi$ un polyn\^ome qui d\'efinit la valuation $\mu$, alors nous d\'eduisons de la proposition \ref{prop:degre-minimal} que tout $b$ appartenant \`a $B$ v\'erifie ${\rm deg}_K(b) \geq {\rm deg}\phi$ et comme $\phi$ a une racine appartenant \`a $B$ nous trouvons ${\rm deg} \phi = {\rm deg}_K(B)$. 

R\'eciproquement soit $b$ appartenant \`a $B$ avec ${\rm deg}_K(b) = {\rm deg}_K(B)$, et soit $\psi$ son polyn\^ome irr\'eductible sur $K$. 
Les polyn\^omes $\phi$ et $\psi$ ont m\^eme degr\'e $d$ et si nous appelons respectivement $d_{\rm s}$ et $d'_{\rm s}$ le nombre de racines distinctes de $\phi$ et $\psi$ nous pouvons \'ecrire $d = p^n d_{\rm s}$ et $d = p^{n'} d'_{\rm s}$. 

Soit $k$ le nombre de boules ferm\'ees $B^{(l)}$ conjugu\'ees \`a $B$, alors d'apr\`es la proposition \ref{prop:valeur-de-mu} nous avons les \'egalit\'es $d_{\rm s} = kc$ et $d'_{\rm s} = kc'$ o\`u $c$ et $c'$ sont respectivement le nombre de racines de $\phi$ et de $\psi$ appartenant \`a une boule $B^{(l)}$. 
Nous en d\'eduisons l'\'egalit\'e 
$$p^nc \ = \ p^{n'}c' \ = \ d/k \ ,$$
et d'apr\`es le corollaire \ref{cor:valeur-de-mu} nous trouvons $\mu(\phi ) = \mu(\psi )$. 
Le th\'eor\`eme est alors une cons\'equence de la proposition \ref{prop:polynome-definissant-la-valuation}.  
\hfill \tf   
\end{preuve}  
 \vskip .2cm

Soit $a$ appartenant \`a $\bar K$, alors pour tout $\delta$ dans $\Gamma _{\bar\nu}$ nous pouvons d\'efinir l'entier $\rho _a(\delta )$ par 
$$\rho _a(\delta ) \ = \ {\rm deg}_K(B(a,\delta )) \ .$$ 
Nous avons ainsi une application croissante $\rho _a$ de $\Gamma _{\bar\nu}$ dans $\N$ major\'ee par $d={\rm deg}_K(a)$.  
De plus la valeur de $\rho _a$ sur $\{ \delta ' \in \Gamma _{\bar\nu} / \delta ' \leq \delta \}$ ne d\'epend pas de $a$ mais uniquement de la boule $B=B(a,\delta )$, en effet si $b$ appartient \`a $B$ alors pour tout $\delta ' \leq \delta$ nous avons encore $B(a,\delta ') = B(b,\delta ')$.  

Pour tout $e\in \N$, avec $e\leq d$ nous d\'efinissons l'ensemble 
$$R_a(e) \ = \ \{ \delta ' \in \Gamma _{\bar\nu} \ / \ \rho _a(\delta ') \leq e \} \ .$$

Nous consid\'erons une valuation bien sp\'ecifi\'ee $\bar\mu$ de $\bar K[x]$ associ\'ee \`a la boule ferm\'ee $B$, sa restriction $\mu$ \`a $K[x]$ et $a$ dans $B$ tel que $\phi$ le polyn\^ome irr\'eductible de $a$ sur $K$ d\'efinisse la valuation $\mu$. 
D'apr\`es ce qui pr\'ec\`ede le degr\'e $d$ de $\phi$ est \'egal \`a ${\rm deg}_K(B)$, et soient $\delta$ et $\gamma$ tels que nous ayons $B=B(a,\delta )$ et $\mu (\phi )=\gamma$.  

\vskip .2cm 

Supposons que l'ensemble $R_a(d-1)$ a un plus grand \'el\'ement $\delta _1$ et soit $d_1=\rho_a(\delta _1)$, alors pour tout $\delta ' \leq \delta _1$ nous avons $\rho_a(\delta ') \leq d_1$ et pour tout $\delta ' > \delta _1$ nous avons $\rho_a(\delta ') =d$. Nous posons $B_1=B(a,\delta _1)$ et nous choisissons $a_1\in B_1$ avec ${\rm deg}_K(a_1) = d_1$. 

Nous appelons $\bar\mu _1$ la valuation de $\bar K[x]$ associ\'ee \`a la boule $B_1$, $\mu _1$ sa restriction \`a $K[x]$ qui est d\'efinie par $\phi _1$ le polyn\^ome irr\'eductible de $a_1$ sur $K$. 

\begin{proposition}\label{prop:mu-valuation-augmentee} 
Le polyn\^ome $\phi$ est un polyn\^ome-cl\'e pour la valuation $\mu _1$ et la valuation $\mu$ est la valuation augment\'ee $[\mu _1 \ ; \ \mu (\phi )=\gamma ]$. 
\end{proposition}

\begin{preuve}
La valuation $\mu _1$ v\'erifie $\mu _1 \leq \mu$ et nous pouvons consid\'erer l'ensemble 
$$\tilde\Phi _{\mu}(\mu _1) \ = \ \{ f \in K[x] \ | \ \mu _1 (f) < \mu (f) \} \ ,$$ 
cet ensemble est non vide car $\phi$ v\'erifie $\mu _1(\phi ) < \mu (\phi )$. En effet si nous \'ecrivons 
$$\phi = \prod _{i=1}^{d_{\rm s}} (x - b_i) ^{p^n} \ ,$$ 
avec $b_1=a$ nous avons $\bar\mu _1(x-b_i) \leq \bar\mu (x-b_i)$ pour tout $i \geq 2$ et $\bar\mu _1(x-a) < \bar\mu (x-a)$. 

Nous appelons $d'$ le degr\'e minimal d'un polyn\^ome appartenant \`a $\tilde\Phi _{\mu}(\mu _1)$ et nous consid\'erons l'ensemble 
$$\Phi _{\mu}(\mu _1) \ = \ \{ \psi \in K[x] \ | \ \mu _1(\psi ) < \mu (\psi ) \ , \ {\rm deg}\psi = d' \ \hbox{et $\psi$ unitaire} \} \ .$$ 

Soit $\psi$ dans $\Phi _{\mu}(\mu _1)$, alors par construction nous avons $d_1 \leq {\rm deg}\psi \leq d$, et d'apr\`es \cite{Va 1}, th\'eor\`eme 1.15, nous savons que $\psi$ est un polyn\^ome-cl\'e pour la valuation $\mu _1$ et que la valuation augment\'ee $\mu ' = [\mu _1 \ ; \ \mu '(\psi )=\gamma ' ]$, avec $\gamma ' = \mu (\psi )$ v\'erifie $\mu _1 \leq \mu ' \leq \mu$. 

Nous d\'eduisons du th\'eor\`eme \ref{th:prolongement-valuation-augmentee} qu'il existe un prolongement $\bar\mu '$ de la valuation $\mu '$ \`a $\bar K[x]$ v\'erifiant $\bar\mu _1 \leq \bar\mu ' \leq \bar\mu$. 
Soit $B'$ la boule associ\'ee \`a $\bar\mu '$, nous avons alors $B \subset B' \subset B_1$, en particulier il existe $\delta '$ tel que $B' = B(a,\delta ')$ et comme la valuation $\mu '$ est diff\'erente de $\mu _1$ nous avons $\delta ' > \delta _1$. 
Par cons\'equent nous avons ${\rm deg}\psi = {\rm deg}\phi$, le polyn\^ome $\phi$ appartient \`a l'ensemble $\Phi _{\mu}(\mu _1)$ et est donc un polyn\^ome-cl\'e pour la valuation $\mu _1$. 

\hfill \tf   
\end{preuve}  
\vskip .2cm

Supposons maintenant que l'ensemble $R_a(d-1)$ n'a pas de plus grand \'el\'ement, il existe alors $\delta _1$ dans $R_a(d-1)$ tel que 
pour tout $\delta '$ appartenant \`a $R_a(d-1)$ avec $\delta ' \geq \delta _1$ nous ayons $\rho_a(\delta ') = \rho_a(\delta _1) = d_1$, et comme pr\'ec\'edemment nous posons $B_1=B(a,\delta _1)$ et nous choisissons $a_1\in B_1$ avec ${\rm deg}_K(a_1) = d_1$. 

Nous choisissons un sous-ensemble $\{\delta _{\alpha} \ / \ \alpha \in A \}$ de $\Gamma _{\bar\nu}$, index\'e par un ensemble totalement ordonn\'e $A$ sans plus grand \'el\'ement, avec $\delta _{\alpha} < \delta _{\alpha '}$ pour $\alpha < \alpha '$ dans $A$, qui soit cofinal dans l'ensemble $\{ \delta ' \in  R_a(d-1) \ / \ \delta ' \geq \delta _1 \}$. 

Nous appelons $\bar\mu _{\alpha}$ la valuation de $\bar K[x]$ associ\'ee \`a la boule $B_{\alpha}= B(a,\delta _{\alpha })$ et $\mu _{\alpha}$ sa restriction \`a $K[x]$. Pour tout $\alpha \in A$ nous choisissons $a_{\alpha}$ dans la boule $B_{\alpha} = B(a,\delta _{\alpha })$ avec ${\rm deg}_K (a_{\alpha })=d_1$, et nous appelons $\phi _{\alpha}$ le polyn\^ome irr\'eductible de $a _{\alpha}$ sur $K$. 

\begin{proposition}\label{prop:mu-valuation-augmentee-limite} 
La famille de valuations ${\mathcal C} = \bigl ( \mu _{\alpha} \bigr ) _{\alpha \in A}$ est une famille continue de valuations de $K[x]$, le polyn\^ome $\phi$ est un polyn\^ome-cl\'e limite pour cette famille et la valuation $\mu$ est la valuation augment\'ee limite $\bigl [ \bigl (  \mu _{\alpha} \bigr ) _{\alpha \in A} \ ; \ \mu (\phi ) = \gamma \bigr ]$. 
\end{proposition}

\begin{preuve}
Tous les polyn\^omes $\phi _{\alpha}$ sont de m\^eme degr\'e et les valuations $\mu _{\alpha}$ v\'erifient $\mu _{\alpha} \leq \mu _{\alpha '}$ pour $\alpha \leq \alpha '$, nous en d\'eduisons que la famille  ${\mathcal C} = \bigl ( \mu _{\alpha} \bigr ) _{\alpha \in A}$ est une famille continue de valuations de $K[x]$. 

De fa\c con analogue \`a la d\'emonstration de la proposition \ref{prop:mu-valuation-augmentee} nous consid\'erons l'ensemble 
$$\tilde\Phi ({\mathcal C}) \ = \ \{ f \in K[x] \ | \ \mu _{\alpha} (f) < \mu _{\alpha '}(f) \ \hbox{pour tout} \ \alpha < \alpha ' \ \hbox{dans} \ A \} \ ,$$ 
qui est non vide car $\phi$ y appartient. Si nous appelons $d_{\mathcal C}$ le degr\'e minimal d'un polyn\^ome dans $\tilde\Phi ({\mathcal C})$, par construction nous avons $d_{\mathcal C} = d ={\rm deg}\phi$, c'est-\`a-dire que $\phi$ appartient \`a l'ensemble 
$$\Phi ({\mathcal C}) \ = \ \{ \psi \in \tilde\Phi ({\mathcal C})\ , \psi \, \hbox{unitaire et deg} \psi = d_{\mathcal C} \} \ ,$$ 
et nous d\'eduisons de la proposition 1.21 de \cite{Va 1} que $\phi$ est un polyn\^ome-cl\'e limite pour la famille ${\mathcal C}$ et que la valuation $\mu$ est la valuation augment\'ee limite $\bigl [ \bigl (  \mu _{\alpha} \bigr ) _{\alpha \in A} \ ; \ \mu (\phi ) = \gamma \bigr ]$. 

\hfill \tf   
\end{preuve}  
\vskip .2cm

\begin{remark}\label{rmq:famille} 
Gr\^ace au corollaire \ref{cor:extensions} nous pouvons associer \`a toute valuation bien sp\'ecifi\'ee $\bar\mu$ de $\bar K[x]$ une famille de valuations ${\mathcal A}_{\bar\mu} = \bigl ( \bar\mu _i \bigr ) _{i \in I}$ v\'erifiant: 
\begin{enumerate}
\item la famille ${\mathcal A} = \bigl ( \mu _i \bigr ) _{i \in I}$ est une famille associ\'ee \`a la valuation $\mu$ de $K[x]$, o\`u nous notons respectivement $\mu$ et $\mu _i$ la valuation de $K[x]$ restriction de la valuation $\bar\mu$ et $\bar\mu _i$;  
\item pour tout $i\leq j$ dans $I$ nous avons $\bar\mu _i \leq \bar\mu _j$, ce qui est \'equivalent \`a $B_j \subset B_i$, o\`u nous appelons $B_i$ la boule ferm\'ee de $\bar K$ associ\'ee \`a la valuation $\bar\mu _i$ de $\bar K[x]$.  
\end{enumerate} 
\end{remark} 
Nous avons construit cette famille \`a partir de la construction de la famille $\mathcal A$ associ\'ee \`a la valuation $\mu$, et la famille ${\mathcal A} = \bigl ( \mu _i \bigr ) _{i \in I}$ est construite en \emph{suivant un ordre croissant}, c'est-\`a-dire que pour d\'eterminer la valuation $\mu _i$ nous avons besoin de conna\^itre les valuations $\mu _j$ pour $j <i$, plus pr\'ecis\'ement nous consid\'erons les ensembles $\tilde\Phi _{\mu}(\mu _j) \ = \ \{ f \in K[x] \ | \ \mu _j (f) < \mu (f) \}$. 

\vskip .2cm 

Nous pouvons utiliser les r\'esultats pr\'ec\'edents pour construire la famille ${\mathcal A}_{\bar\mu} = \bigl ( \bar\mu _i \bigr ) _{i \in I}$ en \emph{suivant un ordre d\'ecroissant}, c'est-\`a-dire en construisant la valuation $\bar\mu _i$ \`a partir des valuations $\mu _j$ pour $j >i$. 

\begin{theorem}\label{th:construction-famille} 
Soit $\bar\mu$ une valuation bien sp\'ecifi\'ee de $\bar K[x]$, alors il existe une famille ${\mathcal A}_{\bar\mu} = \bigl ( \bar\mu _i \bigr ) _{i \in I}$ de valuations de $\bar K[x]$ v\'erifiant les propri\'et\'es (1) et (2) de la remarque \ref{rmq:famille} obtenue de la fa\c con suivante. 

D'apr\`es la propri\'et\'e (2) les boules $B_i$ associ\'ees aux valuations $\bar\mu _i$ sont de la forme $B_i=B(a, \delta _i)$ avec une famille $\bigl ( \delta _i \bigr ) _{i\in I}$ d\'ecroissante. 
Si nous connaissons la valuation $\bar\mu _i$ de la famille, les valuations $\bar\mu _j$ pour $j<i$ sont d\'etermin\'ees en consid\'erant l'ensemble 
$$R_K(\bar\mu _i) \ = \ \bigl \{ \delta \ / \ {\rm deg}_K (B(a,\delta )) < {\rm deg}_K (B(a,\delta _i)) \bigr \} \ .$$    

Si cet ensemble a un plus grand \'el\'ement $\delta _{i-1}$ nous construisons la valuation $\bar\mu _{i-1}$ comme la valuation associ\'ee \`a la boule $B(a,\delta _{i-1})$, et $i-1$ est le pr\'ed\'ecesseur de $i$ dans $I$.  

Si cet ensemble n'a pas de plus grand \'el\'ement nous choisissons une famille $(\delta _{\alpha} )_{\alpha \in A}$ cofinale dans $R_K(\bar\mu _i)$ et nous d\'efinissons la famille $(\bar\mu _{i_{\alpha}})_{i_{\alpha} \in I_A}$ o\`u la valuation $\bar\mu _{i_{\alpha}}$ est la valuation associ\'ee \`a la boule $B(a,\delta _{\alpha})$, et $I_A=( i_{\alpha}) _{\alpha \in A}$ est la sous-famille de $I$ sans plus grand \'el\'ement cofinale dans $\{ j\in I / j<i \}$. 

Nous nous arr\^etons, c'est-\`a-dire nous trouvons la premi\`ere valuation $\bar\mu _1$ de la famille quand nous trouvons un ensemble $R_K(\bar\mu _i)$ tel que ${\rm deg}_K (B(a,\delta ))=1$, o\`u $\delta = \delta _{i-1}$ ou $\delta = \delta _{\alpha}$, avec les notations pr\'ec\'edentes. 
\end{theorem} 

\begin{preuve}
C'est une cons\'equence imm\'ediate des propositions \ref{prop:mu-valuation-augmentee}  et \ref{prop:mu-valuation-augmentee-limite}.  

\hfill \tf   
\end{preuve}  
\vskip .2cm

\begin{appendix}

\section{Suites pseudo-convergentes et extension imm\'ediate}\label{appA}

Dans cette partie nous allons montrer comment il est possible d'interpr\'eter les r\'esultats de Kaplansky sur les familles pseudo-convergentes et les extensions imm\'ediates (\cite{Ka}) \`a partir des notions de familles admissibles continues et de valuations augment\'ees limites. 

Soit ${\mathcal C} = \bigl ( \mu _{\alpha} \bigr ) _{\alpha \in A}$ une famille continue de valuations, nous notons respectivement $\bigl ( \phi _{\alpha} \bigr ) _{\alpha \in A}$ et $\bigl ( \gamma _{\alpha} \bigr ) _{\alpha \in A}$ les familles de polyn\^omes-cl\'es et de valeurs dans $\Gamma$ associ\'ees. 
Nous pouvons supposer que l'ensemble $A$ a un plus petit \'el\'ement $\vartheta _0$, nous notons $\mu$ la valuation $\mu _{\vartheta _0}$ et toute valuation $\mu _{\alpha}$ pour $\alpha > \vartheta _0$ est la valuation augment\'ee $\mu _{\alpha} = [ \mu \ ; \ \mu _{\alpha} (\phi _{\alpha})=\gamma _{\alpha} ]$, et nous notons $d$ le degr\'e des polyn\^omes $\phi _{\alpha}$. 
Rappelons le r\'esultat suivant qui d\'ecrit le comportement des valeurs $\mu _{\alpha} (f)$ pour un polyn\^ome $f$ de $K[x]$ quand $\alpha$ parcourt l'ensemble $A$, et qui est une cons\'equence directe du th\'eor\`eme 1.19 de \cite{Va 1}. 

\begin{theorem} \label{th:comportement-des-mu-de-f} 
Soit $f$ dans $K[x]$ de degr\'e $m$, il existe un entier $v = v_{\mathcal C}(f)$ v\'erifiant $0 \leq v \leq k= [m/d]$, une famille $\bigl ( \alpha _i \bigr ) _{0\leq i \leq j}$ dans $A$ avec $\vartheta _0 = \alpha _0 \leq \alpha _1 \leq \ldots \leq \alpha _j$ et $\delta$ dans $\Gamma$, o\`u $j=k-v$, tels que 
$$\mu _{\alpha} (f)\  = \ \delta + \sum _{i=1} ^{j} Inf(\gamma _{\alpha} , \gamma _{\alpha _i}) + v \gamma _{\alpha} \ .$$
\end{theorem} 

En particulier $f$ appartient \`a l'ensemble $\tilde\Phi \bigl ( {\mathcal C} \bigr ) =  \{ f \in K[x] \ | \ \mu _{\alpha} (f) < \mu _{\beta}(f) , \forall \alpha < \beta \in A \}$ si et seulement si $v_{\mathcal C}(f) > 0$. 

\vskip .2cm 

\begin{corollary}\label{cor:comportement-des-mu-de-f} 
La fonction 
$$
\begin{array}{cccc}
\mu _.(f) : & A & \longrightarrow & \Gamma \\
 & \alpha & \mapsto & \mu _{\alpha} (f)
\end{array}$$ 
est lin\'eaire par morceaux et concave. 

Plus pr\'ecis\'ement nous avons 
$$\begin{array} {rl}
 \mu _{\alpha} (f)  \ = \ \delta _i + (k+1-i) \gamma _{\alpha} & \hbox{pour}\quad  \alpha _{i-1} \leq \alpha \leq \alpha _i \ , 1\leq i \leq j\ , \\ 
 \mu _{\alpha} (f)  \ = \ \delta _{j +1}+ (k-j) \gamma _{\alpha} & \hbox{pour}\quad   \alpha \geq \alpha _j\ , 
\end{array}$$ 
o\`u $\delta _i = \delta + \sum _{l=1} ^{i-1} \gamma _l$. 
\end{corollary} 

\vskip .2cm

Soit $f$ dans $K[x]$ et pour tout $\alpha$ nous consid\'erons la division euclidienne de $f$ par $\phi _{\alpha}$ que nous notons $f = q_{\alpha} \phi _{\alpha} + r_{\alpha}$. 
Comme le polyn\^ome $\phi _{\beta}$ est un polyn\^ome-cl\'e pour la valuation $\mu _{\alpha}$ pour tout $\beta \geq \alpha$, d'apr\`es le lemme 1.1 de \cite{Va 1}, nous avons l'in\'egalit\'e: 
$$\mu (r _{\beta}) = \mu _{\alpha} (r _{\beta}) \geq \mu _{\alpha} (f) \ ,$$ 
avec l'in\'egalit\'e stricte si et seulement si $f$ est $\mu _{\alpha}$-divisible par $\phi _{\beta}$.   
\vskip .2cm

\begin{theorem}\label{th:reste} 
Soit $f$ un polyn\^ome dans $K[x]$ et soit $f = q_{\alpha} \phi _{\alpha} + r_{\alpha}$ la division euclidienne de $f$ par le polyn\^ome-cl\'e $\phi _{\alpha}$: 

\begin{enumerate} 
\item si $f$ n'appartient pas \`a $\tilde\Phi \bigl ( {\mathcal C} \bigr )$, il existe $\alpha _1$ tel que pour tout $\alpha > \alpha _1$ nous avons 
$$\mu _{\alpha} (f) = \mu ( r_{\alpha}) \ .$$

\item si $f$ appartient \`a $\tilde\Phi \bigl ( {\mathcal C} \bigr )$, il existe $\alpha _1$ tel que pour tout $\beta > \alpha > \alpha _1$ nous avons 
$$\mu _{\alpha} (f) \ \leq \ \mu (r _{\alpha}) \ < \ \mu _{\beta} (f)  \ . $$
 \end{enumerate}  
\end{theorem}

\vskip .2cm

\begin{corollary}\label{cor:reste}
\begin{enumerate}
\item Si $f$ n'appartient pas \`a $\tilde\Phi \bigl ( {\mathcal C} \bigr )$, il existe $\alpha _1$ tel que pour tout $\beta > \alpha > \alpha _1$ nous avons $\mu (r_{\beta}) = \mu ( r_{\alpha})$. 

\item Si $f$ appartient \`a $\tilde\Phi \bigl ( {\mathcal C} \bigr )$, il existe $\alpha _1$ tel que pour tout $\beta > \alpha > \alpha _1$ nous avons $\mu (r_{\beta}) > \mu ( r_{\alpha})$. 
\end{enumerate}
\end{corollary} 

\vskip .2cm

\begin{preuveth}  
Si $f$ n'appartient pas \`a $\tilde\Phi \bigl ( {\mathcal C} \bigr )$, il existe $\alpha _1$ tel que pour tout $\beta \geq \alpha \geq \alpha _1$ nous avons l'\'egalit\'e $\mu _{\beta}(f) = \mu _{\alpha}(f)$, et de plus le polyn\^ome $f$ n'est pas $\mu _{\alpha}$-divisible par $\phi _{\beta}$, d'o\`u l'\'egalit\'e $\mu (r _{\beta}) = \mu _{\alpha} (f)$. 

\vskip .2cm 

Si $f$ appartient \`a $\tilde\Phi \bigl ( {\mathcal C} \bigr )$, nous d\'eduisons du corollaire \ref{cor:comportement-des-mu-de-f} qu'il existe $\alpha _{[f]}$ dans $A$, $\delta _{[f]}$ dans $\Gamma$ et un entier $v _{[f]} \geq 1$ tels que pour tout $\beta \geq \alpha _{[f]}$ nous ayons l'\'egalit\'e $\mu _{\beta} (f) = \delta _{[f]} + v _{[f]} \gamma _{\beta}$. 

Pour tout $\beta \leq \alpha$, le polyn\^ome $\phi _{\alpha}$ est un polyn\^ome-cl\'e pour la valuation $\mu _{\beta}$, par cons\'equent nous avons $\mu _{\beta} (f) \leq \mu ( r_{\alpha})$, de plus pour $\beta < \alpha$ nous avons l'in\'egalit\'e $\mu _{\beta} (f) < \mu _{\alpha} (f)$, par cons\'equent $f$ est $\mu _{\beta}$-divisible par $\phi _{\alpha}$ d'o\`u l'in\'egalit\'e $\mu (r _{\alpha}) > \mu _{\beta} (f)$.  
Nous en d\'eduisons que pour tout $\beta$ avec $\alpha _{[f]} < \beta < \alpha$ nous avons $\mu _{\beta} (q _{\alpha}) = \delta _{[f]} + ( v_{[f]} -1) \gamma _{\beta}$. 
Nous d\'eduisons du corollaire \ref{cor:comportement-des-mu-de-f} que nous avons $\mu _{\beta} (q _{\alpha}) \leq \delta _{[f]} + ( v_{[f]} -1) \gamma _{\beta}$ pour tout $\beta > \alpha _{[f]}$, par cons\'equent pour tout $\beta > \alpha$ nous avons l'in\'egalit\'e stricte 
$$\mu _{\beta} ( q _{\alpha} \phi _{\alpha} ) \ \leq \ \delta _{[f]} + ( v_{[f]} -1) \gamma _{\beta} + \gamma _{\alpha} \ < \ \delta _{[f]} + v_{[f]} \gamma _{\beta} \ = \ \mu _{\beta} (f) \ ,$$
et nous en d\'eduisons $\mu _{\beta} (f) > \mu ( r_{\alpha})$. 

\hfill \tf    
\end{preuveth}   

\vskip .2cm

\begin{corollary}\label{cor:reste}
Si le groupe des valeurs $\Gamma$ ne poss\`ede de pas de plus petit \'el\'ement strictement positif, c'est-\`a-dire si le sous-groupe isol\'e minimal de $\Gamma$ n'est pas discret, il existe $\alpha _1$ tel que pour tout $\alpha \geq \alpha _1$ nous avons 
$$\mu _{\alpha}(f) = \mu ( r _{\alpha}) \ . $$
\end{corollary} 

\begin{preuve}
Gr\^ace au lemme 1.17 de \cite{Va 1} nous pouvons choisir la famille $\mathcal C$ de telle fa\c con que toute valeur de $\Gamma$ sup\'erieure \`a $\mu _{\alpha}(f)$ soit atteinte par $\mu _{\beta}(f)$, pour un $\beta$ dans $A$. 
Par cons\'equent si $\Gamma$ ne poss\`ede de pas de plus petit \'el\'ement strictement positif nous d\'eduisons le r\'esultat des in\'egalit\'es $\mu _{\beta} (f) > \mu ( r_{\alpha}) \geq \mu _{\alpha} (f)$. 

\hfill \tf    
\end{preuve}   

\vskip .2cm

Nous avons vu qu'il est \'equivalent de se donner une famille pseudo-convergente $( a _{\alpha }) _{\alpha \in A}$ d'\'el\'ements de $K$ et de se donner une famille continue ${\mathcal C} = \bigl ( \mu _{\alpha} \bigr ) _{\alpha \in A}$ de valuations de $K[x]$ telle que les polyn\^omes $\phi _{\alpha}$ qui la d\'efinissent sont de degr\'e un, c'est-\`a-dire sont de la forme $\phi _{\alpha} = x - a_{\alpha}$ (cf. proposition \ref{prop:corpsalgclos}).
Dans ce cas pour tout polyn\^ome $f$ de $K[x]$ le reste de la division de $f$ par le polyn\^ome-cl\'e $\phi _{\alpha} = x - a _{\alpha}$ est \'egal \`a $f(a_{\alpha})$. 

Nous d\'eduisons alors de ce qui pr\'ec\`ede le r\'esultat suivant. 

\begin{corollary}\label{cor:reformulation}
Soit $( a _{\alpha }) _{\alpha \in A}$ une famille pseudo-convergente d'\'el\'ements de $K$, et soit $f$ un polyn\^ome de $K[x]$, alors  il existe $\alpha _1$ tel que pour tout $\beta > \alpha > \alpha _1$ nous avons:  
\begin{enumerate}
\item si $f$ n'appartient pas \`a $\tilde\Phi \bigl ( {\mathcal C} \bigr )$, alors $\nu (f (a_{\beta})) = \nu (f ( a_{\alpha}))$, 
\item si $f$ appartient \`a $\tilde\Phi \bigl ( {\mathcal C} \bigr )$, alors $\nu (f (a_{\beta})) > \nu (f ( a_{\alpha}))$. 
\end{enumerate}
\end{corollary} 

\vskip .2cm 

Nous consid\'erons maintenant un corps valu\'e $(K,\nu )$ et nous voulons \'etudier ses extensions imm\'ediates, c'est-\`a-dire les extensions de corps valu\'es $(L,\nu _L ) / (K,\nu )$ telles que les extensions de groupes  $\Gamma _{\nu _L} / \Gamma _{\nu}$ et les extensions r\'esiduelles $\kappa _{\nu _L} / \kappa _{\nu}$ soient triviales.

Nous remarquons d'abord que d'apr\`es la proposition \ref{prop:egalite-d-Abhyankar} si $\mu$ est une valuation bien sp\'ecifi\'ee de $K(x)$ v\'erifiant $\Gamma _{\nu} = \Gamma _{\mu}$ l'extension r\'esiduelle $\kappa _{\mu} / \kappa _{\nu}$ est de degr\'e de transcendance un, en particulier n'est pas triviale. Par cons\'equent les extensions monog\`enes imm\'ediates $(L,\nu _L )$ de $(K,\nu )$ sont d\'efinies soit par des pseudo-valuations, c'est le cas d'une extension alg\'ebrique $L$ de $K$, soit par une valuation de $K[x]$ qui n'est pas bien sp\'ecifi\'ee, c'est le cas $L=K(x)$. 

\begin{proposition}\label{prop:extension-immediate}  
Si $(L,\nu _L )$ est une extension imm\'ediate monog\`ene de $(K,\nu )$, la famille admise ${\mathcal A}$ associ\'ee \`a la valuation ou pseudo-valuation $\mu$ de $K[x]$ d\'efinie par $\nu _L$ est de la forme suivante:
\begin{enumerate} 
\item si $L$ est une extension transcendante pure, $L=K(x)$, la famille ${\mathcal A}$ est une famille admissible continue ${\mathcal C}=\bigl ( \mu _{\alpha } \bigr ) _{\alpha \in A}$ telle que les polyn\^omes $\phi _{\alpha}$ sont de degr\'e un et telle que l'ensemble $\tilde\Phi \bigl ( {\mathcal C} \bigr )$ est vide; 

\item si $L$ est une extension alg\'ebrique, $L=K(a )\simeq K[x]/(\phi )$ avec $\phi$ polyn\^ome irr\'eductible de $a$ sur $K$, la famille ${\mathcal A}$ est compos\'ee d'une famille admissible continue ${\mathcal C}=\bigl ( \mu _{\alpha } \bigr ) _{\alpha \in A}$ telle que les polyn\^omes $\phi _{\alpha}$ sont de degr\'e un et de la pseudo-valuation $\mu$, o\`u $\mu$ est obtenue comme valuation augment\'ee-limite de la famille ${\mathcal C}$ associ\'ee au polyn\^ome-cl\'e limite $\phi$. 
\end{enumerate}  
\end{proposition}  

\begin{preuve}  
Si la famille $\mathcal A$ associ\'ee \`a la valuation, ou pseudo-valuation $\mu$ contient un couple de valuations successives $(\mu _{i-1}, \mu _i )$, c'est-\`a-dire telle que la valuation $\mu _i$ est obtenue comme valuation augment\'ee $\mu _i = [ \mu _{i-1} \ ; \ \mu _i(\phi _i )=\gamma _i ]$, nous d\'eduisons de la proposition 2.5 et de la d\'emonstration de la proposition 2.9 de \cite{Va 4} que nous avons l'in\'egalit\'e 
$$e(\nu _L/\nu ) f(\nu _L/\nu ) \geq {\rm deg}\phi _i/{\rm deg}\phi_{i-1}\ ,$$
o\`u nous notons respectivement $e(\nu _L/\nu )$ et $f(\nu _L/\nu )$ l'indice de ramification et le degr\'e de l'extension r\'esiduelle de $(L,\nu _L= /(K,\nu )$, et o\`u $\phi _{i-1}$ et $\phi _i$ sont les polyn\^omes d\'efinissant les valuations $\mu _{i-1}$ et $\mu _i$. 

En particulier si l'extension $(L,\nu _L) /(K,\nu )$ est imm\'ediate nous en d\'eduisons que pour tout couple de valuations successives $(\mu _{i-1}, \mu _i )$ de la famille $\mathcal A$ nous avons l'\'egalit\'e ${\rm deg}\phi _i = {\rm deg}\phi_{i-1}$, par cons\'equent que la famille ne contient pas de partie discr\`ete  . 

Dans le cas o\`u la valuation $\mu$ n'est pas bien sp\'ecifi\'ee, la famille admise $\mathcal A$ associ\'ee \`a $\mu$ est ouverte, elle est constitu\'ee d'une seule famille simple ${\mathcal S}$ qui est de la forme ${\mathcal S} = \bigl ( ( \mu _{\alpha})_{\alpha \in A}\bigr )$, avec $\tilde\Phi \bigl ( ( \mu _{\alpha}) _{\alpha \in A} \bigr ) = \emptyset$ (cf. remarque \ref{rmq:famille ouverte}), et de plus nous d\'eduisons de ce qui pr\'ec\`ede que tous les polyn\^omes-cl\'es $\phi _{\alpha}$ sont de degr\'e un. 

Dans le cas o\`u la valuation $\mu$ est bien sp\'ecifi\'ee, $\mu$ est une pseudo-valuation obtenue comme valuation augment\'ee limite de la famille continue ${\mathcal S} = \bigl ( ( \mu _{\alpha})_{\alpha \in A}\bigr )$, telle que les polyn\^omes-cl\'es $\phi _{\alpha}$ sont de degr\'e un, $\mu = \bigl [ ( \mu _{\alpha})_{\alpha \in A} \ ; \ \mu (\phi ) = +\infty \bigr ]$. 

\hfill \tf    
\end{preuve}

\vskip .2cm

\begin{proposition}\label{prop:ensemble-Phi} 
Soit $(L,\nu _L)$ une extension imm\'ediate monog\`ene de $(K,\nu )$ et soit ${\mathcal C} =\bigl ( \mu _{\alpha} \bigr ) _{\alpha \in A}$ la famille admissible continue associ\'ee d\'efinie \`a la proposition \ref{prop:extension-immediate}. 
Alors un polyn\^ome $f$ n'appartient pas \`a l'ensemble $\tilde\Phi \bigl ( {\mathcal C}\bigr )$ si et seulement si il existe $\alpha$ dans $A$ tel que nous ayons $\nu _L(f) = \nu (f(a_{\beta}))$ pour tout $\beta \geq \alpha$.   
\end{proposition}  

\begin{preuve}  
Si le polyn\^ome $f$ n'appartient pas \`a l'ensemble $\tilde\Phi \bigl ( {\mathcal C} \bigr )$ il existe $\alpha$ dans $A$ tel que pour tout $\beta \geq \alpha$ nous ayons $\mu _{\alpha}(f) = \mu _{\beta} (f) = \nu _L(f)$, c'est-\`a-dire tels que $f$ n'est pas $\mu _{\alpha}$-divisible par $\phi _{\beta}$. 
Par cons\'equent d'apr\`es le corollaire \ref{cor:reformulation}, nous avons pour tout $\beta\geq \alpha$ l'\'egalit\'e $\nu _L(f) = \nu (f(a_{\beta}))$. 

\hfill \tf    
\end{preuve}

Nous pouvons d\'eduire de ce qui pr\'ec\`ede les r\'esultats suivants, qui sont des reformulations des th\'eor\`emes 2 et 3 de \cite{Ka}.

\begin{corollary} 
Soit $(L,\nu _L)$ une extension imm\'ediate monog\`ene de $(K,\nu )$ et soit $\bigl ( a _{\alpha} \bigr ) _{\alpha \in A}$ la famille pseudo-convergente associ\'ee, alors l'extension $L/K$ est transcendante si et seulement si pour tout $f$ dans $K[x]$ il existe $\alpha$ dans $A$ tel que nous ayons $\nu _L(f) = \nu (f(a_{\beta}))$ pour tout $\beta \geq \alpha$.  
\end{corollary} 

\begin{preuve}  
Nous sommes dans le cas d'une extension transcendante si et seulement si l'ensemble $\tilde\Phi \bigr ( {\mathcal C}\bigr )$ est vide, c'est-\`a-dire si et seulement si pour tout $f$ dans $K[x]$ il existe $\alpha$ dans $A$ tel que nous ayons $\mu _{\beta}(f) = \nu (f(a_{\beta}))$ pour tout $\beta \geq \alpha$. 

\hfill \tf    
\end{preuve}

\begin{corollary} 
Soit $(L,\nu _L)$ une extension imm\'ediate monog\`ene de $(K,\nu )$ et soit $\bigl ( a _{\alpha} \bigr ) _{\alpha \in A}$ la famille pseudo-convergente associ\'ee, et nous supposons que l'extension $L/K$ est alg\'ebrique. 
Alors si $\phi$ est un polyn\^ome v\'erifiant $\nu (f(a_{\beta})) > \nu (f(a_{\alpha}))$ de degr\'e minimal, l'extension $L$ est l'extension $L=K(z)$ o\`u $z$ est une racine du polyn\^ome $\phi$. 

De plus pour tout $y$ dans $L$ il existe $r(x)$ dans $K[x]$ avec ${\rm deg}r < {\rm deg}\phi$   tel que $y=r(z)$ et la valuation $\nu _L$ est d\'efinie par $\nu _L(y) = \nu ( r(a_{\alpha}))$ pour $\alpha$ assez grand. 
\end{corollary} 

\begin{preuve}  
Nous sommes dans le cas d'une extension alg\'ebrique si et seulement si l'ensemble $\tilde\Phi \bigr ( {\mathcal C}\bigr )$ est non vide, tout polyn\^ome unitaire $\phi$ de degr\'e minimal dans $\tilde\Phi \bigr ( {\mathcal C}\bigr )$ est un polyn\^ome-cl\'e limite pour la famille continue ${\mathcal C} = \bigl ( \mu _{\alpha} \bigr ) _{\alpha \in A}$  et la pseudo-valuation augment\'ee limite $\mu$ de $K[x]$ d\'efinie par $\mu = \bigl [ \bigl ( \mu _{\alpha} \bigr ) _{\alpha \in A} \ ; \ \mu (\phi ) = + \infty \bigr ]$  induit une valuation $\nu _L$ de l'extension $L=K[x] / ( \phi )$ telle que l'extension est imm\'ediate. 

De plus tout \'el\'ement $y$ de $L$ est d\'efini par un polyn\^ome $r$ de degr\'e strictement inf\'erieur au degr\'e de $\phi$, par cons\'equent $r$ n'appartient pas \`a $\tilde\Phi \bigr ( {\mathcal C}\bigr )$ et nous trouvons $\nu _L(y) = \mu (r) = \mu _{\alpha} (r)$ pour $\alpha$ assez grand. 

\hfill \tf    
\end{preuve}

\end{appendix}

\vskip 2cm 

		  

\vskip 1cm

\end{document}